\documentclass[12pt]{amsart}
\usepackage{amsmath}
\usepackage{amssymb}
\usepackage{amsfonts}
\usepackage{amsthm}
\usepackage{amscd}
\usepackage{bm}
\usepackage[dvips]{graphicx}
\newtheorem{thrm}{Theorem}

\newtheorem{thm}{Theorem}[section]
\newtheorem{lem}[thm]{Lemma}
\newtheorem{pro}[thm]{Proposition}

\newtheorem{prob}[thm]{Problem}
\newcommand{\C}{\mathbf{C}}

\newcommand{\Z}{\mathbf{Z}}
\newcommand{\Q}{\mathbf{Q}}
\newcommand{\N}{\mathbf{N}}

\begin{document}
\title[On the graded quotients of the $\mathrm{SL}(m,\C)$-representation algebras of groups]
{On the graded quotients of the $\mathrm{SL}(m,\C)$-representation algebras of groups}
\address[Takao Satoh]{Department of Mathematics, Faculty of Science Division II, Tokyo University of Science,
         1-3 Kagurazaka, Shinjuku, Tokyo, 162-8601, Japan}
\email{takao@rs.tus.ac.jp}
\maketitle

\begin{center}
{\sc Takao Satoh}\footnote{e-address: takao@rs.tus.ac.jp}
\end{center}

\begin{abstract}
In this paper, we consider certain descending filtrations of the $\mathrm{SL}(m,\C)$-representation algebras of free groups and
free abelian groups. By using it, we introduce analogs of the Johnson homomorphisms of the automorphism groups of free groups.
We show that the first homomorphisms are extended to the automorphism groups of free groups as crossed homomorphisms.
Furthermore we show that the extended crossed homomorphisms induce Kawazumi's cocycles and Morita's
cocycles.
This works are generalization of our previous results for the $\mathrm{SL}(2,\C)$-representation algebras.
\end{abstract}
\section{Introduction}\label{S-Int}

For any $m \geq 2$ and any group $G$, let $R^m(G)$ be the set $\mathrm{Hom}(G, \mathrm{SL}(m,\C))$ of all $\mathrm{SL}(m,\C)$-representations of $G$.
Let $\mathcal{F}(R^m(G),\C)$ be the set $\{ \chi : R^m(G) \rightarrow \C \}$ of all complex-valued functions on $R^m(G)$.
Then $\mathcal{F}(R^m(G),\C)$ naturally has the $\C$-algebra structure coming from the pointwise sum and product. (See Subsection {\rmfamily \ref{S-Rep_R}} for details.)
For any $x \in G$ and any $1 \leq i, j \leq m$, we define the element $a_{ij}(x)$ of $\mathcal{F}(R^m(G),\C)$ to be
\[ (a_{ij}(x))(\rho) := {\text{the $(i,j)$-component of $\rho(x)$}} \]
for any $\rho \in R^m(G)$. We call the map $a_{ij}(x)$ the $(i,j)$-component function of $x$, or simply a component function of $x$.
Let $\mathfrak{R}_{\Q}^m(G)$ be the $\Q$-subalgebra of $\mathcal{F}(R^m(G),\C)$ generated by $a_{ij}(x)$ for all $x \in G$ and $1 \leq i, j \leq m$.
We call $\mathfrak{R}_{\Q}^m(G)$ the $\mathrm{SL}(m,\C)$-representation algebras of $G$ over $\Q$.
In this paper, we introduce a descending filtration of $\mathfrak{R}_{\Q}^m(G)$ consisting of $\mathrm{Aut}\,G$-invariant ideals, and study the
graded quotients of it.

\vspace{0.5em}

To the best of our knowledge, the study of the algebra $\mathfrak{R}_{\Q}^m(G)$ has a not so long history.
Let $F_n$ be a free group of rank $n$ with basis $x_1, \ldots, x_n$.
Classically, the $\Q$-subalgebra of $\mathfrak{R}_{\Q}^2(F_n)$ generated by characters of $F_n$ was actively studied.
For any $x \in F_n$, the map $\mathrm{tr}\,x := a_{11}(x) + a_{22}(x)$ is called the Fricke character of $x$.
Fricke and Klein \cite{Fri} used the Fricke characters to study the classification of Riemann surfaces.
In the 1970s, Horowitz \cite{Ho1} and \cite{Ho2} investigated several algebraic properties of the algebra of Fricke characters
by using the combinatorial group theory.
In 1980, Magnus \cite{Mg2} studied some relations among Fricke characters of free groups systematically.
As is found in Acuna-Maria-Montesinos's paper \cite{AJM}, today Magnus's research has been developed to the study of the $\mathrm{SL}(2,\C)$-character
varieties of free groups by quite many authors.
Let $\mathfrak{X}_{\Q}^2(F_n)$ be the $\Q$-subalgebra of $\mathfrak{R}_{\Q}^2(F_n)$ generated by $\mathrm{tr}\,x$ for all $x \in F_n$.
The algebra $\mathfrak{X}_{\Q}^2(F_n)$ is called the algebra of Fricke characters of $F_n$.
Let $\mathfrak{C}$ be the ideal of $\mathfrak{X}_{\Q}^2(F_n)$ generated by $\mathrm{tr}\,x - 2$ for all $x \in F_n$.
In our previous papers \cite{HS1} and \cite{S19}, we considered an application of the theory of the Johnson homomorphisms of $\mathrm{Aut}\,F_n$
by using the Fricke characters. In particular, we determined the structure of the graded quotients
$\mathrm{gr}^k(\mathfrak{C}) := \mathfrak{C}^k/\mathfrak{C}^{k+1}$ for $1 \leq k \leq 2$, introduced analogs of the Johnson homomorphisms,
and showed that the first homomorphism extends to $\mathrm{Aut}\,F_n$ as a crossed homomorphism.

\vspace{0.5em}

We briefly review the history of the Johnson homomorphisms.
In 1965, Andreadakis \cite{And} introduced a certain descending central filtration of $\mathrm{Aut}\,F_n$ by using the natural action of
$\mathrm{Aut}\,F_n$ on the nilpotent quotients of $F_n$.
We call this filtration the Andreadakis-Johnson filtration of $\mathrm{Aut}\,F_n$.
(For definitions, see Section {\rmfamily \ref{Ss-A_J}}.)
In the 1980s, Johnson studied such filtration for mapping class groups of surfaces in order to investigate the
group structure of the Torelli groups in a series of works \cite{Jo1}, \cite{Jo2}, \cite{Jo3} and \cite{Jo4}.
In particular, he determined the abelianization of the Torelli group by introducing a certain homomorphism.
Today, his homomorphism is called the first Johnson homomorphism, and it
is generalized to higher degrees.
Over the last two decades, the Johnson homomorphisms of the mapping class groups have been actively studied from various viewpoints
by many authors including Morita \cite{Mo1}, Hain \cite{Hai} and others.

\vspace{0.5em}

The Johnson homomorphisms are naturally defined for $\mathrm{Aut}\,F_n$.
\[ \tilde{\tau}_k : \mathcal{A}_{F_n}(k) \hookrightarrow \mathrm{Hom}_{\Z}(H, \mathcal{L}_{F_n}(k+1)) \]
where $H := F_n/[F_n,F_n]$ is the abelianization of $F_n$.
(See Subsection {\rmfamily \ref{Ss-A_J}} for details.)
So far, we concentrate on the study of the cokernels of Johnson homomorphisms of $\mathrm{Aut}\,F_n$ in a series of our works \cite{S03}, \cite{S06}, \cite{S11} and \cite{ES1}
with combinatorial group theory and representation theory.
Since each of the Johnson homomorphisms is a $\mathrm{GL}(n,\Z)$-equivariant injective homomorphism,
it is an important problem to determine its images and cokernels. 
On the other hand, the study of the extendability of the Johnson homomorphisms have been received attentions.
Morita \cite{Mo5} showed that the first Johnson homomorphism of the mapping class group, which initial domain is the Torelli group,
can be extended to the mapping class group as a crossed homomorphism
by using the extension theory of groups.
Inspired by Morita's work, Kawazumi \cite{Kaw} obtained a corresponding results for $\mathrm{Aut}\,F_n$
by using the Magnus expansion of $F_n$. Furthermore, he constructed higher twisted cohomology classes with the extended first Johnson homomorphism and the cup product.
By restricting them to the mapping class group, he investigated relations between the higher cocycles and the Morita-Mumford classes.
Recently, Day \cite{Day} showed that each of Johnson homomorphisms of $\mathrm{Aut}\,F_n$ can be extended to a crossed homomorphism from $\mathrm{Aut}\,F_n$
into a certain finitely generated free abelian group.

\vspace{0.5em}

As mentioned above, we \cite{HS1} constructed analogs of the Johnson homomorphisms with the algebra $\mathfrak{X}_{\Q}^2(F_n)$ of Fricke characters,
and showed that the first homomorphism can be extended to $\mathrm{Aut}\,F_n$ in \cite{S19}.
It is, however, difficult to push forward with our research since
the structures of the graded quotients $\mathrm{gr}^k(\mathfrak{C})$ are too complicated to handle.
In \cite{S20}, we considered a similar situation for the $\mathrm{SL}(2,\C)$-representation algebra $\mathfrak{R}_{\Q}^2(F_n)$.
In this paper, we generalize our previous works to the $\mathrm{SL}(m,\C)$-representation case.
Set $s_{ij}(x) := a_{ij}(x) - \delta_{ij}$
for any $1 \leq i, j \leq m$ and $x \in F_n$ where $\delta$ means Kronecker's delta.
Let $\mathfrak{J}_{F_n}$ be the ideal of $\mathfrak{R}_{\Q}^m(F_n)$ generated by $s_{ij}(x_l)$ for all $1 \leq i, j \leq m$ and $1 \leq l \leq n$.
Then the products of $\mathfrak{J}_{F_n}$ define a descending filtration of $\mathfrak{R}_{\Q}^m(F_n)$: 
$\mathfrak{J}_{F_n} \supset \mathfrak{J}_{F_n}^2 \supset \mathfrak{J}_{F_n}^3 \supset \cdots$, 
which consists of $\mathrm{Aut}\,F_n$-invariant ideals. Set $\mathrm{gr}^k(\mathfrak{J}_{F_n}) := \mathfrak{J}_{F_n}^k/\mathfrak{J}_{F_n}^{k+1}$ for any $k \geq 1$. The first aim of the paper is to give a basis of $\mathrm{gr}^k(\mathfrak{J}_{F_n})$ as a $\Q$-vector space.
Set
\[ T_k := \Big{\{} \prod_{\substack{ 1 \leq i, j \leq m \\[1pt] (i,j) \neq (m,m)}} \prod_{l=1}^n s_{ij}(x_l)^{e_{ij,l}}
    \,\Big{|}\, e_{ij,l} \geq 0, \,\,\,
    \sum_{\substack{ 1 \leq i, j \leq m \\[1pt] (i,j) \neq (m,m)}} \sum_{l=1}^n e_{ij,l} = k \Big{\}}
   \subset \mathfrak{J}_{F_n}^k. \] 
\begin{thrm}[Propositions {\rmfamily \ref{P-Rythm}} and {\rmfamily \ref{P-dim-2}}]\label{I-T-1} \,\
For each $k \geq 1$, the set $T_k \pmod{\mathfrak{J}_{F_n}^{k+1}}$ forms a basis of $\mathrm{gr}^k(\mathfrak{J}_{F_n})$ as a $\Q$-vector space.
Furthermore, for any $n \geq 2$ and $k \geq 1$, we have
\[ \mathrm{gr}^k(\mathfrak{J}_{F_n}) \cong {\bigoplus}^{'}
    \bigotimes_{\substack{ 1 \leq i, j \leq m \\[1pt] (i,j) \neq (m,m)}} S^{e_{ij}} H_{\Q} \]
as a $\mathrm{GL}(n,\Z)$-module. Here the sum runs over all tuples $(e_{ij})$ for $1 \leq i, j \leq m$ and $(i,j) \neq (m,m)$ such that
the sum of the $e_{ij}$ is equal to $k$.
\end{thrm}
\noindent
This theorem is a generalization of our previous result for the case where $k=2$ in \cite{S20}.

\vspace{0.5em}

Now, set
\[ \mathcal{D}_{F_n}^m(k) := \mathrm{Ker}(\mathrm{Aut}\,F_n \rightarrow \mathrm{Aut}(\mathfrak{J}_{F_n}/\mathfrak{J}_{F_n}^{k+1})). \]
The groups $\mathcal{D}_{F_n}^m(k)$ define a descending central filtration of $\mathrm{Aut}\,F_n$.
Let $\mathcal{A}_{F_n}(1) \supset \mathcal{A}_{F_n}(2) \supset \cdots$ be the Andreadakis-Johnson filtration of $\mathrm{Aut}\,F_n$.
The second purpose of the paper is to show a relation between $\mathcal{A}_{F_n}(k)$ and $\mathcal{D}_{F_n}^m(k)$ as follows.
\begin{thrm}[Theorems {\rmfamily \ref{T-Marine}} and {\rmfamily \ref{T-Flower}}]\label{I-T-2} \,\
\begin{enumerate}
\item[$(1)$] For any $k \geq 1$, $\mathcal{A}_{F_n}(k) \subset \mathcal{D}_{F_n}^m(k)$.
\item[$(2)$] For any $k \geq 1$ and $m \geq 2$, we have $\mathcal{D}_{F_n}^{m+1}(k) \subset \mathcal{D}_{F_n}^{m}(k)$.
\end{enumerate}
\end{thrm}
In \cite{S20}, we showed that $\mathcal{D}_{F_n}^2(k)=\mathcal{A}_{F_n}(k)$ for $1 \leq k \leq 4$.
Hence, we have $\mathcal{D}_{F_n}^m(k)=\mathcal{A}_{F_n}(k)$ for any $1 \leq k \leq 4$.
By referring to the theory of the Johnson homomorphisms, we can construct analogs of them:
\[ \tilde{\eta}_k : \mathcal{D}_{F_n}^m(k) \rightarrow \mathrm{Hom}_{\Q}(\mathrm{gr}^1(\mathfrak{J}_{F_n}), \mathrm{gr}^{k+1}(\mathfrak{J}_{F_n})) \]
by $\tilde{\eta}_k(\sigma) :=(f \pmod{\mathfrak{J}_{F_n}^2} \mapsto f^{\sigma} - f \pmod{\mathfrak{J}_{F_n}^{k+2}}$)
for any $\sigma \in \mathcal{D}_{F_n}^m(k)$.
(See Subsection {\rmfamily \ref{Ss-DFD}} for details.)
In this paper, we consider an extension of the first homomorphism $\tilde{\eta}_1$,
and study some relations to the extension of the first Johnson homomorphism $\tilde{\tau}_1$ of $\mathrm{Aut}\,F_n$. 

\vspace{0.5em}

Set $H_{\Q} := H \otimes_{\Z} \Q$.
In \cite{S01}, we computed $H^1(\mathrm{Aut}\,F_n, H_{\Q}) = \Q$, and showed that it is generated by Morita's cocycle
$f_M$. On the other hand, we \cite{S10} also computed
$H^1(\mathrm{Aut}\,F_n, H_{\Q}^* \otimes_{\Q} \Lambda^2 H_{\Q}) =\Q^{\oplus 2}$, and showed that
it is generated by Kawazumi's cocycle $f_K$ and the cocycle induced from $f_M$. Kawazumi's cocycle $f_K$ is an extension
of $\tilde{\tau}_1$.
(See Subsection {\rmfamily \ref{Ss-1-coh}} for details.)
In Subsection {\rmfamily \ref{Ss-1-coh}}, we construct the crossed homomorphism
\[ \theta_{F_n} : \mathrm{Aut}\,F_n \rightarrow \mathrm{Hom}_{\Q}(\mathrm{gr}^1(\mathfrak{J}_{F_n}), \mathrm{gr}^2(\mathfrak{J}_{F_n})) \]
which is an extension of $\tilde{\eta}_1$. By taking suitable reductions of the target of $\theta_{F_n}$,
we obtain the crossed homomorphisms
\[ f_1 : \mathrm{Aut}\,F_n \rightarrow H_{\Q}^* \otimes_{\Q} \Lambda^2 H_{\Q}, \hspace{1em}
   f_2 : \mathrm{Aut}\,F_n \rightarrow H_{\Q}. \]
The third purpose of the paper is to show the following theorem.
\begin{thrm}[Theorem {\rmfamily \ref{T-Haruhi}}]\label{I-T-3}
For any $n \geq 2$,
\[\begin{split}
 f_K = f_1, \hspace{1em} f_M = - f_2 + \delta_x
  \end{split}\]
for $x=x_1+x_2 + \cdots + x_n \in H_{\Q}$ as crossed homomorphisms. Here $\delta_x$ is the principle derivation associated with $x$.
\end{thrm}
\noindent
This shows that our crossed homomorphism induces both of Kawazumi's cocycle and Morita's cocycle, and that
$\theta_{F_n}$ defines the non-trivial cohomology class in
$H^1(\mathrm{Aut}\,F_n, $ $\mathrm{Hom}_{\Q}(\mathrm{gr}^1(\mathfrak{J}_{F_n}), \mathrm{gr}^2(\mathfrak{J}_{F_n})))$.

\vspace{0.5em}

In \cite{HS2} and \cite{S22}, we studied $\mathfrak{X}_{\Q}^2(H)$ and $\mathfrak{R}_{\Q}^2(H)$.
In this paper, we generalize the results in \cite{S22} to the $\mathrm{SL}(m,\C)$-representation cases.
By using a parallel argument, we can define a descending filtration
$\mathfrak{J}_{H} \supset \mathfrak{J}_{H}^2 \supset \mathfrak{J}_{H}^3 \supset \cdots$
of ideals in $\mathfrak{R}_{\Q}^m(H)$.
In contrast with the free group case, however, it is quite hard to determine the structure of the graded quotients $\mathrm{gr}^k(\mathfrak{J}_H)$.
Here, we gave basis of $\mathrm{gr}^k(\mathfrak{J}_H)$ for $1 \leq k \leq 2$. In particular, we see
\[ \mathrm{gr}^1(\mathfrak{J}_H) \cong H_{\Q}^{\oplus m^2-1},
 \hspace{1em} \mathrm{gr}^2(\mathfrak{J}_H) \cong (S^2 H_{\Q})^{\oplus \frac{1}{2}m^2(m^2-1)} \oplus (\Lambda^2 H_{\Q})^{\{ \oplus \frac{1}{2}(m^2-1)(m^2-4) \} }. \]
We remark that for a general $m \geq 3$, the situation of the $\mathrm{SL}(m,\C)$-representation case
is much more different and complicated than those of the $\mathrm{SL}(2,\C)$-representation case.
At the present stage, we have no idea to give a result for a general $k \geq 3$.
By using the above results, we construct the crossed homomorphism
\[ \theta_H : \mathrm{Aut}\,F_n \rightarrow \mathrm{Hom}_{\Q}(\mathrm{gr}^1(\mathfrak{J}_H), \mathrm{gr}^2(\mathfrak{J}_H)), \]
and show that it induces Morita's cocycle $f_M$ in Theorem {\rmfamily \ref{T-Haruhi-H}}.

\vspace{0.5em}

In Section {\rmfamily \ref{S-Rep_R}},
we study the $\mathrm{SL}(m,\C)$-representation algebra for any group $G$, and introduce the descending central filtration
$\mathcal{D}_G^m(1) \supset \mathcal{D}_G^m(2) \supset \cdots$ of $\mathrm{Aut}\,G$.
In Section {\rmfamily \ref{S-Ext-J}}, we discuss a general construction of a crossed homomorphism of $\mathrm{Aut}\,G$
by using an associative algebra on which $\mathrm{Aut}\,G$ acts.
Then, in Sections {\rmfamily \ref{S-Free}} and {\rmfamily \ref{S-FAb}}, we apply the above results for the free group case
and the free abelian group case respectively.

\tableofcontents

\section{Notation and conventions}\label{S-Not}

\vspace{0.5em}

Throughout the paper, we use the following notation and conventions. 
\begin{itemize}
\item Let $G$ be a group. The automorphism group $\mathrm{Aut}\,G$ of $G$ acts on $G$ from the right.
      For any $\sigma \in \mathrm{Aut}\,G$ and $x \in G$, the action of $\sigma$ on $x$ is denoted by $x^{\sigma}$.
\item Let $N$ be a normal subgroup of a group $G$. For an element $g \in G$, we also denote the coset class of $g$ by $g \in G/N$
      if there is no confusion.
      Similarly, for an algebra $R$, an element $f \in R$ and an ideal $I$ of $R$, we also denote by $f$ the coset class of $f$ in $R/I$ if
      there is no confusion.
\item For elements $x$ and $y$ in $G$, the commutator bracket $[x,y]$ of $x$ and $y$
      is defined to be $xyx^{-1}y^{-1}$.
\item On the direct product $\N \times \N$, we consider the usual lexicographic ordering among tuples $(i,j)$ for $i, j \in \N$. The ordering is denoted by
      $\leq_{\mathrm{lex}}$.
\end{itemize}

\vspace{0.5em}

\section{The $\mathrm{SL}(m,\C)$-representation algebras of groups}\label{S-Rep_R}

\vspace{0.5em}

Let $G$ be a group generated by $x_1, \ldots, x_n$. For any $m \geq 2$, we denote by
$R^m(G)$ the set $\mathrm{Hom}(G, \mathrm{SL}(m,\C))$ of all $\mathrm{SL}(m,\C)$-representations of $G$. Let
$\mathcal{F}(R^m(G),\C)$
be the set $\{ \chi : R^m(G) \rightarrow \C \}$ of all complex-valued functions on $R^m(G)$.
Then $\mathcal{F}(R^m(G),\C)$ has the $\C$-algebra structure by the pointwise sum and product defined by
\[\begin{split}
  (\chi+ \chi')(\rho) & := \chi(\rho) + \chi'(\rho), \\
  (\chi \chi')(\rho) & := \chi(\rho) \chi'(\rho), \\
  (\lambda \chi)(\rho) & := \lambda (\chi(\rho)) 
  \end{split}\]
for any $\chi$, $\chi' \in \mathcal{F}(R^m(G),\C)$, $\lambda \in \C$, and $\rho \in R^m(G)$.
The automorphism group $\mathrm{Aut}\,G$ of $G$ naturally acts on $R^m(G)$ and $\mathcal{F}(R^m(G),\C)$ from the right by
\[ \rho^{\sigma}(x) := \rho(x^{\sigma^{-1}}), \hspace{1em} \rho \in R^m(G) \,\,\, \text{and} \,\,\, x \in G \]
and
\[ \chi^{\sigma}(\rho) := \chi(\rho^{\sigma^{-1}}), \hspace{1em} \chi \in \mathcal{F}(R^m(G),\C) \,\,\, \text{and} \,\,\, \rho \in R^m(G) \]
for any $\sigma \in \mathrm{Aut}\,G$.

\vspace{0.5em}

For any $x \in G$ and any $1 \leq i, j \leq m$, we define the element $a_{ij}(x)$ of $\mathcal{F}(R^m(G),\C)$ to be
\[ (a_{ij}(x))(\rho) := {\text{the $(i,j)$-component of $\rho(x)$}} \]
for any $\rho \in R^m(G)$.
The action of an element $\sigma \in \mathrm{Aut}\,G$ on $a_{ij}(x)$ is given by $a_{ij}(x^{\sigma})$.
We have the following relations:
\begin{equation}
a_{ij}(x^{-1}) = \tilde{a}_{ji}(x) \label{eq-Lumy}
\end{equation}
and
\begin{equation}
 a_{ij}(xy) = \sum_{k=1}^m a_{ik}(x)a_{kj}(y) \label{eq-Lumy2}
\end{equation}
for any $1 \leq i, j \leq m$ and $x, y \in G$. Here $\tilde{a}_{ij}(x)$ denotes the $(i,j)$-cofactor of the matrix $(a_{ij}(x))$,
and is written as a polynomial of $a_{kl}(x)$ since the determinant of $(a_{ij}(x))$ is equal to one.

\vspace{0.5em}

Let $\mathfrak{R}_{\Q}^m(G)$ be the $\Q$-subalgebra of $\mathcal{F}(R^m(G),\C)$ generated by $a_{ij}(x)$ for all $x \in G$ and $1 \leq i, j \leq m$.
We call $\mathfrak{R}_{\Q}^m(G)$ the $\mathrm{SL}(m,\C)$-representation algebras of $G$ over $\Q$.
\begin{lem}\label{L-Black}
The algebra $\mathfrak{R}_{\Q}^m(G)$ is finitely generated by $a_{ij}(x_l)$ for $1 \leq l \leq n$.
\end{lem}
\textit{Proof.}
For any $x \in G$, we have $x= x_{i_1}^{e_1} x_{i_2}^{e_2} \cdots x_{i_r}^{e_r}$ for some $1 \leq i_m \leq n$ and $e_j= \pm 1$.
We show $a_{ij}(x) \in \mathfrak{R}_{\Q}^m(G)$ by the induction on $r \geq 1$.
For $r=1$, it is obvious that $a_{ij}(x)=a_{ij}(x_{i_1}) \in \mathfrak{R}_{\Q}^m(G)$ for $e_1=1$, and that
$a_{ij}(x^{-1})=\tilde{a}_{ji}(x_{i_1}) \in \mathfrak{R}_{\Q}^m(G)$ for $e_1=-1$ from (\ref{eq-Lumy}).
For $r \geq 2$, from (\ref{eq-Lumy2}) we have
\[ a_{ij}(x) = \sum_{k=1}^m a_{ik}(x_{i_1}^{e_1} x_{i_2}^{e_2} \cdots x_{i_{r-1}}^{e_{r-1}})a_{kj}(x_{i_r}^{e_r}) \]
Hence by the inductive hypothesis we see $a_{ij}(x) \in \mathfrak{R}_{\Q}^m(G)$. $\square$

\vspace{0.5em}

Next, we consider an $\mathrm{Aut}\,G$-invariant ideal of $\mathfrak{R}_{\Q}^m(G)$.
Set $s_{ij}(x) := a_{ij}(x) - \delta_{ij}$
for any $1 \leq i, j \leq m$ and $x \in G$ where $\delta$ means the Kronecker's delta.
Let $\mathfrak{J}_G$ be the ideal of $\mathfrak{R}_{\Q}^m(G)$ generated by $s_{ij}(x_l)$ for all $1 \leq i, j \leq m$ and $1 \leq l \leq n$.
\begin{lem}\label{L-White}
The ideal $\mathfrak{J}_G$ is $\mathrm{Aut}\,G$-invariant.
\end{lem}
\textit{Proof.}
Since $s_{ij}(x_l)^{\sigma} = s_{ij}(x_l^{\sigma})$ for any $\sigma \in \mathrm{Aut}\,G$,
it suffices to show that $s_{ij}(x) \in \mathfrak{J}_G$ for any $x \in G$.

\vspace{0.5em}

From Lemma {\rmfamily \ref{L-Black}}, $s_{ij}(x)$ is written as a polynomial among $a_{ij}(x_l)$.
Then by substituting $s_{ij}(x) + \delta_{ij}$ for $a_{ij}(x)$, we see that
\[ s_{ij}(x) = C + \text{(terms of degree $\geq 1$)} \]
for some $C \in \Q$.
Let ${\bf 1} : G \rightarrow \mathrm{SL}(m,\C)$ be the trivial representation. 
By observing the image of ${\bf 1}$ by $s_{ij}(x)$, we obtain $C=0$ since $s_{ij}(x_l)({\bf 1})=0$ for any $1 \leq l \leq n$.
This means $s_{ij}(x) \in \mathfrak{J}_G$. Therefore, $\mathfrak{J}_G$ is $\mathrm{Aut}\,G$-invariant. $\square$

\vspace{0.5em}

Then we have the descending filtration
\[ \mathfrak{J}_G \supset \mathfrak{J}_G^2 \supset \mathfrak{J}_G^3 \supset \cdots \]
which consists of $\mathrm{Aut}\,G$-invariant ideals. Set $\mathrm{gr}^k(\mathfrak{J}_G) := \mathfrak{J}_G^k/\mathfrak{J}_G^{k+1}$ for any $k \geq 1$.
Each graded quotient $\mathrm{gr}^k(\mathfrak{J}_G)$ is an $\mathrm{Aut}\,G$-invariant finite dimensional $\Q$-vector space.
Furthermore, the graded sum
\[ \mathrm{gr}(\mathfrak{J}_G) := \bigoplus_{k \geq 1} \mathrm{gr}^k(\mathfrak{J}_G) \]
naturally has the graded $\Q$-algebra structure given by the product map
\[ \mathrm{gr}^k(\mathfrak{J}_G) \times \mathrm{gr}^l(\mathfrak{J}_G) \rightarrow \mathrm{gr}^{k+l}(\mathfrak{J}_G) \]
defined by
\[ ( f \pmod{\mathfrak{J}_G^{k+1}}, \,\, g \pmod{\mathfrak{J}_G^{l+1}}) \mapsto fg \pmod{\mathfrak{J}_G^{k+l+1}}. \]
Recall that
\begin{equation}\label{eq-pooh}
\begin{split}
 s_{ij}(xy) & = s_{ij}(x) + s_{ij}(y) + \sum_{k=1}^m s_{ik}(x)s_{kj}(y) \\
  & \equiv s_{ij}(x) + s_{ij}(y) \pmod{\mathfrak{J}_G^2}, \\
 s_{ij}(x^{-1}) & = - s_{ij}(x) + \sum_{k=1}^m s_{ik}(x)s_{kj}(x^{-1}) \\
  & \equiv - s_{ij}(x) \pmod{\mathfrak{J}_G^2}
\end{split}
\end{equation}
For any $x, y \in G$. 

\vspace{0.5em}

In order to describe generators of $\mathrm{gr}^k(\mathfrak{J}_G)$ as a $\Q$-vector space,
we consider relations among $a_{ij}(x_l)$. Set
\begin{equation}
 \Delta_l := \det{(a_{ij}(x_l))} = \sum_{\sigma \in \mathfrak{S}_m} \mathrm{sgn}(\sigma) a_{1 \sigma(1)}(x_l) a_{2 \sigma(2)}(x_l) \cdots a_{m \sigma(m)}(x_l) \in \mathfrak{R}_{\Q}^m(G) \label{eq-3}
\end{equation}
for any $1 \leq l \leq n$. Since we consider representations into $\mathrm{SL}(m,\C)$, we see that $\Delta_l = 1$ for any $1 \leq l \leq n$.
Thus,
\[\begin{split}
 \Delta_l -1 &= -1 + (s_{11}(x_l)+1) (s_{22}(x_l)+1) \cdots (s_{mm}(x_l)+1) \\
  & \hspace{2em} + \sum_{\substack{\sigma \in \mathfrak{S}_m \\[1pt] \sigma \neq 1}} \mathrm{sgn}(\sigma)
    (s_{1 \sigma(1)}(x_l) + \delta_{1 \sigma(1)}) \cdots(s_{m \sigma(m)}(x_l) + \delta_{m \sigma(m)}) \\
  & = s_{11}(x_l) + s_{22}(x_l) + \cdots + s_{mm}(x_l) + \text{(terms of degree $\geq 2$)} \in \mathfrak{J}_G, 
\end{split}\]
and hence
\begin{equation}
 s_{11}(x_l) + s_{22}(x_l) + \cdots + s_{mm}(x_l) \equiv 0 \pmod{\mathfrak{J}_G^2} \label{eq-bloom}
\end{equation}
The first purpose of the paper is to give a generating set of $\mathrm{gr}^k(\mathfrak{J}_G)$ as a $\Q$-vector space.
For any $k \geq 1$, set
\[ T_k := \Big{\{} \prod_{\substack{ 1 \leq i, j \leq m \\[1pt] (i,j) \neq (m,m)}} \prod_{l=1}^n s_{ij}(x_l)^{e_{ij,l}}
    \,\Big{|}\, e_{ij,l} \geq 0, \,\,\,
    \sum_{\substack{ 1 \leq i, j \leq m \\[1pt] (i,j) \neq (m,m)}} \sum_{l=1}^n e_{ij,l} = k \Big{\}}
   \subset \mathfrak{J}_G^k. \] 
\begin{lem}\label{L-Egret}
For each $k \geq 1$, $T_k$ generates $\mathrm{gr}^k(\mathfrak{J}_G)$ as a $\Q$-vector space.
\end{lem}
\textit{Proof.}
Since $\mathrm{gr}^k(\mathfrak{J}_G)$ is generated by
\[ \prod_{1 \leq i, j \leq m} \prod_{l=1}^n s_{ij}(x_l)^{e_{ij,l}} \hspace{0.5em} \text{for} \hspace{0.5em}
    \sum_{1 \leq i, j \leq m} \sum_{l=1}^n e_{ij,l} = k, \]
by using (\ref{eq-bloom}) to eliminate $s_{mm}(x_l)$, we see that $T_k$ generates $\mathrm{gr}^k(\mathfrak{J}_G)$. $\square$

\begin{lem}\label{L-Dream}
Let $\varphi : G \rightarrow G'$ be a group homomorphism. Then $\varphi$ induces the algebra homomorphism
\[ \overline{\varphi} : \mathfrak{R}_{\Q}^m(G) \rightarrow \mathfrak{R}_{\Q}^m(G') \]
defined by
\[ a_{ij}(x) \mapsto a_{ij}(\varphi(x)) \]
for any $x \in G$. Furthermore, if $\varphi$ is surjective, so is $\overline{\varphi}$.
\end{lem}
\textit{Proof.}
For a given homomorphism $\varphi : G \rightarrow G'$, $\varphi$ induces the map
$\widetilde{\varphi} : R^m(G') \rightarrow R^m(G)$ defined by $\rho \mapsto \rho \circ \varphi$. Then $\widetilde{\varphi}$ induces the algebra homomorphism
$\overline{\varphi} : \mathcal{F}(R^m(G),\C) \rightarrow \mathcal{F}(R^m(G'),\C)$ defined by
\[ \overline{\varphi}(\chi) := \chi \circ \widetilde{\varphi} \]
for any $\chi \in R^m(G')$. The restriction of $\overline{\varphi}$ to $\mathfrak{R}_{\Q}^m(G)$ is the required homomorphism.
$\square$

\section{Descending filtrations of $\mathrm{Aut}\,G$}\label{S-DFL}

\subsection{The Andreadakis-Johnson filtration}\label{Ss-A_J}
\hspace*{\fill}\ 

\vspace{0.5em}

In this subsection, we review the Andreadakis-Johnson filtration of $\mathrm{Aut}\,G$ without proofs. The main purpose of the subsection is to fix the notations.
For basic materials concerning the Andreadakis-Johnson filtration and the Johnson homomorphisms of $\mathrm{Aut}\,G$, see \cite{S06} or \cite{S15}, for example.

\vspace{0.5em}

First, consider the lower central series $\Gamma_G(1) \supset \Gamma_G(2) \supset \cdots$ of $G$ defined by the rule
\[ \Gamma_G(1):= G, \hspace{1em} \Gamma_G(k) := [\Gamma_G(k-1),G], \hspace{1em} k \geq 2. \]
For any $y_1, \ldots, y_k \in G$, the left-normed commutator
\[ [[ \cdots [[y_1, y_2], y_3], \cdots ], y_k] \]
of weight $k$ is denoted by
\[ [y_1, y_2, \cdots, y_k ] \]
for simplicity, and called a simple $k$-fold commutator. Then we have
\begin{lem}[See Section 5.3 in \cite{MKS}.]\label{L-Mint}
For any $k \geq 1$, the group $\Gamma_G(k)$ is generated by all simple $k$-fold commutators.
\end{lem}
\noindent
For any $k \geq 1$, set $\mathcal{L}_G(k) := \Gamma_G(k)/\Gamma_G(k+1)$.
\begin{lem}\label{L-Pine}
If $G$ is generated by $x_1, \ldots , x_n$, then each of the graded quotients
$\mathcal{L}_G(k)$ is generated by (the coset class of) the simple $k$-fold commutators
\[ [x_{i_1},x_{i_2}, \ldots , x_{i_k}], \hspace{1em} 1 \leq i_j \leq n \]
as an abelian group.
\end{lem}
\noindent
For a proof, see Theorem 5.4 in \cite{MKS}, for example. This shows that if $G$ is finitely generated then so is $\mathcal{L}_G(k)$
for any $k \geq 1$.

\vspace{0.5em}

For $k \geq 1$, the action of $\mathrm{Aut}\,G$ on each nilpotent quotient $G/\Gamma_G(k+1)$ induces the homomorphism
$\mathrm{Aut}\,G \rightarrow \mathrm{Aut}(G/\Gamma_G(k+1))$.
We denote its kernel by $\mathcal{A}_G(k)$. Then the groups $\mathcal{A}_G(k)$ define the descending filtration
\[ \mathcal{A}_G(1) \supset \mathcal{A}_G(2) \supset \cdots \supset \mathcal{A}_G(k) \supset \cdots \]
of $\mathrm{Aut}\,G$. We call this the Andreadakis-Johnson filtration of $\mathrm{Aut}\,G$.
The first term $\mathcal{A}_G(1)$ is called the IA-automorphism group of $G$, and is also denoted by $\mathrm{IA}(G)$.
Namely, $\mathrm{IA}(G)$ consists of automorphisms which act on the abelianization $G^{\mathrm{ab}}$ of $G$ trivially.

\vspace{0.5em}

The Andreadakis-Johnson filtration of $\mathrm{Aut}\,G$ was originally introduced by Andreadakis \cite{And} in the 1960s.
The name \lq\lq Johnson" comes from Dennis Johnson who studied the Johnson filtration and the Johnson homomorphism for the
mapping class group of a surface in the 1980s.
In particular, Andreadakis showed that
\begin{thm}[Andreadakis, \cite{And}]\label{T-And} \quad
\begin{enumerate}
\item For any $k$, $l \geq 1$, $\sigma \in \mathcal{A}_G(k)$ and $x \in \Gamma_G(l)$, $x^{-1} x^{\sigma} \in \Gamma_G(k+l)$.
\item For any $k$ and $l \geq 1$, $[\mathcal{A}_G(k), \mathcal{A}_G(l)] \subset \mathcal{A}_G(k+l)$.
\item If $\displaystyle \bigcap_{k \geq 1} \Gamma_G(k) =1$, then $\displaystyle \bigcap_{k \geq 1} \mathcal{A}_G(k) =1$.
\end{enumerate}
\end{thm}

\vspace{0.5em}

To study the structure of the graded quotients ${\mathrm{gr}}^k(\mathcal{A}_G) :=\mathcal{A}_G(k)/\mathcal{A}_G(k+l)$,
we use the Johnson homomorphisms of $\mathrm{Aut}\,G$.
For any $\sigma \in \mathcal{A}_G(k)$, consider the map
$\tilde{\tau}_k(\sigma) :G^{\mathrm{ab}} \rightarrow \mathcal{L}_G(k+1)$ defined by
\[ x \pmod{\Gamma_G(2)} \hspace{1em} \mapsto \hspace{1em} x^{-1}x^{\sigma} \pmod{\Gamma_G(k+2)} \]
for any $x \in G$. Then $\tilde{\tau}_k(\sigma)$ is a homomorphism between abelian groups.
Furthermore, the map $\tilde{\tau}_k : \mathcal{A}_G(k) \rightarrow \mathrm{Hom}_{\Z}(G^{\mathrm{ab}}, \mathcal{L}_G(k+1))$ defined by
\[ \sigma \hspace{0.5em} \mapsto \hspace{0.5em} \tilde{\tau}_k(\sigma) \]
is a homomorphism. From the definition, it is easy to see that the kernel of $\tilde{\tau}_k$ is $\mathcal{A}_G(k+1)$.
Therefore $\tilde{\tau}_k$ induces the injective homomorphism
\[ \tau_k : \mathrm{gr}^k (\mathcal{A}_G) \hookrightarrow \mathrm{Hom}_{\Z}(G^{\mathrm{ab}}, \mathcal{L}_G(k+1)). \]
For each $k \geq 1$, we call both of the homomorphisms $\tilde{\tau}_k$ and $\tau_k$ the $k$-th Johnson homomorphisms of $\mathrm{Aut}\,G$.

\subsection{The filtration $\mathcal{D}_G^m(1) \supset \mathcal{D}_G^m(2) \supset \cdots$}\label{Ss-DFD}
\hspace*{\fill}\ 

\vspace{0.5em}

For $k \geq 1$, the action of $\mathrm{Aut}\,G$ on each nilpotent quotient $\mathfrak{J}_G/\mathfrak{J}_G^{k+1}$ induces the homomorphism
$\mathrm{Aut}\,G \rightarrow \mathrm{Aut}(\mathfrak{J}_G/\mathfrak{J}_G^{k+1})$.
Set
\[ \mathcal{D}_G^m(k) := \mathrm{Ker}(\mathrm{Aut}\,G \rightarrow \mathrm{Aut}(\mathfrak{J}_G/\mathfrak{J}_G^{k+1})). \]
The groups $\mathcal{D}_G^m(k)$ define the descending filtration
\[ \mathcal{D}_G^m(1) \supset \mathcal{D}_G^m(2) \supset \cdots \supset \mathcal{D}_G^m(k) \supset \cdots \]
of $\mathrm{Aut}\,G$.
For any $f \in \mathfrak{J}_G$ and $\sigma \in \mathrm{Aut}\,G$, set
\[ s_{\sigma}(f) := f^{\sigma}- f \in \mathfrak{J}_G. \]
In the following, we give a few lemmas and a proposition without proofs. 
In \cite{S20}, we give proofs of the lemmas and the proposition for the case where $m=2$ and $G=F_n$.
We show the followings by the same way as in \cite{S20}. For details, see \cite{S20}.
\begin{lem}\label{L-Rouge}
For any $f \in \mathfrak{J}_G$ and $\sigma, \tau \in \mathrm{Aut}\,G$,
\begin{enumerate}
\item $s_{\sigma \tau}(f) = (s_{\sigma}(f))^{\tau} + s_{\tau}(f)$,
\item $s_{1_{G}}(f) = 0$,
\item $s_{\sigma^{-1}}(f) = - (s_{\sigma}(f))^{\sigma^{-1}}$,
\item $s_{[\sigma, \tau]}(f) = \{ s_{\tau}(s_{\sigma}(f)) - s_{\sigma}(s_{\tau}(f)) \}^{\sigma^{-1}\tau^{-1}}$.
\end{enumerate}
\end{lem}

\begin{lem}\label{L-Peach}
For any $k, l \geq 1$, $f \in \mathfrak{J}_G^l$ and $\sigma \in \mathcal{D}_G^m(k)$, we have $s_{\sigma}(f) \in \mathfrak{J}_G^{k+l}$.
\end{lem}

\begin{pro}\label{P-Passion}
For any $k, l \geq 1$, $[\mathcal{D}_G^m(k), \mathcal{D}_G^m(l)] \subset \mathcal{D}_G^m(k+l)$.
\end{pro}
\noindent
This proposition shows that the filtration $\mathcal{D}_G^m(1) \supset \mathcal{D}_G^m(2) \supset \cdots$ is a central filtration of $\mathrm{Aut}\,G$.
Next, we show that $\mathcal{D}_G^m(k)$ contains $\mathcal{A}_G(k)$ for any $k \geq 1$.

\begin{lem}\label{L-Pine}
For any $1 \leq  i, j \leq m$, $k \geq 1$ and $y \in \Gamma_G(k)$, we have $s_{ij}(y) \in \mathfrak{J}_G^k$.
\end{lem}
\textit{Proof.}
From the fact that $\Gamma_G(k)$ is generated by all simple $k$-fold commutators, and from (\ref{eq-pooh}),
it suffices to show the lemma in the case where $y=[y_1, y_2, \ldots, y_k]$
for any $y_1, \ldots, y_k \in G$.
We use the induction on $k \geq 1$. If $k=1$, the lemma is obvious. Assume $k \geq 2$ and set $z:=[y_1, \ldots, y_{k-1}]$.
Then $s_{ij}(z) \in \mathfrak{J}_G^{k-1}$ by the inductive hypothesis. Furthermore, from (\ref{eq-Lumy}), we see
$s_{ij}(z^{-1}) \in \mathfrak{J}_G^{k-1}$.
By using (\ref{eq-pooh}), we have
\[\begin{split}
   s_{ij}([z,y_k]) & = s_{ij}(z y_k z^{-1} y_k^{-1}) \\
     & = s_{ij}(zy_k) + s_{ij}(z^{-1} y_k^{-1}) + \sum_{h=1}^m s_{ih}(zy_k) s_{hj}(z^{-1} y_k^{-1}) \\
     & \equiv s_{ij}(z) + s_{ij}(y_k) + s_{ij}(z^{-1}) + s_{ij}(y_k^{-1}) + \sum_{h=1}^m s_{ih}(y_k) s_{hj}(y_k^{-1}) \pmod{\mathfrak{J}_G^k} \\
     & \equiv s_{ij}(z) + s_{ij}(z^{-1}) \pmod{\mathfrak{J}_G^k} \\
     & \equiv 0 \pmod{J^k}.
  \end{split}\]
Here, the last equation follows from
\[ 0 = s_{ij}(1_G) = s_{ij}(zz^{-1}) \equiv s_{ij}(z) + s_{ij}(z^{-1}) \pmod{\mathfrak{J}_G^k}. \]
$\square$

\begin{lem}\label{L-Berry}
For any $k \geq 1$, $z \in G$ and $y \in \Gamma_G(k)$, we have $s_{ij}(zy) \equiv s_{ij}(y) \pmod{\mathfrak{J}_G^k}$.
\end{lem}
\textit{Proof.}
This lemma follows from (\ref{eq-pooh}) and Lemma {\rmfamily \ref{L-Pine}} immediately. $\square$

\vspace{0.5em}

From Lemma {\rmfamily \ref{L-Berry}}, we obtain the following theorem.
\begin{thm}\label{T-Marine}
For any $k \geq 1$, $\mathcal{A}_G(k) \subset \mathcal{D}_G^m(k)$.
\end{thm}

\vspace{0.5em}

Now, we give a relation between $\mathcal{D}_G^m(k)$ and $\mathcal{D}_G^{m+1}(k)$ for $m \geq 2$.
\begin{thm}\label{T-Flower}
For any $k \geq 1$ and $m \geq 2$, we have $\mathcal{D}_G^{m+1}(k) \subset \mathcal{D}_G^{m}(k)$.
\end{thm}
\textit{Proof.}
Take any $\sigma \in \mathcal{D}_G^{m+1}(k)$ and $x \in G$. For any $1 \leq i, j \leq m$, in order to distinguish $s_{ij}(x) \in \mathfrak{R}_{\Q}^{m+1}(G)$
from $s_{ij}(x) \in \mathfrak{R}_{\Q}^{m}(G)$, we write $s_{ij}^m(x)$ for $s_{ij}(x)$ if we consider $s_{ij}(x)$ as an element in $\mathfrak{R}_{\Q}^{m}(G)$.
For any representation $\rho : G \rightarrow \mathrm{SL}(m,\C)$, define the representation $\overline{\rho} : G \rightarrow \mathrm{SL}(m+1,\C)$ by
\[ \overline{\rho}(x) := \begin{pmatrix} \rho(x) & 0 \\ 0 & 1 \end{pmatrix}. \]
Then we see
\[\begin{split}
   s_{ij}^m(x^{\sigma})(\rho) &= s_{ij}^m(\rho(x^{\sigma})) = s_{ij}^{m+1}(\overline{\rho}(x^{\sigma})) = s_{ij}^{m+1}(x^{\sigma})(\overline{\rho}) \\
      & = (s_{ij}^{m+1}(x) +( \text{polynomial among $s_{pq}^{m+1}(y)$ of degree $\geq k+1$}))(\overline{\rho}) \\
      & = s_{ij}^{m+1}(\overline{\rho}(x)) + (\text{polynomial among $s_{pq}^{m+1}(\overline{\rho}(y))$ of degree $\geq k+1$}) \\
      & = s_{ij}^m(\rho(x)) + (\text{polynomial among $s_{pq}^{m}(\rho(y))$ of degree $\geq k+1$}) \\
      & = s_{ij}^m(x)(\rho) + (\text{polynomial among $s_{pq}^{m}(\rho(y))$ of degree $\geq k+1$})
  \end{split}\]
for some $y \in G$.
Hence we see $s_{ij}^m(x^{\sigma}) \equiv s_{ij}^m(x) \pmod{\mathfrak{J}_G^{k+1}}$ in $\mathfrak{R}_{\Q}^{m}(G)$. This shows $\sigma \in \mathcal{D}_G^{m}(k)$.
$\square$

\vspace{0.5em}

From the above theorem, we have
\[ \mathcal{A}_G(k) \subset \cdots \subset \mathcal{D}_G^{m+1}(k) \subset \mathcal{D}_G^{m}(k) \subset \cdots \subset \mathcal{D}_G^{2}(k). \]
Here we give a sufficient condition for $\mathcal{A}_G(k)=\mathcal{D}_G^m(k)$. For any $1 \leq i, j \leq m$ and $k \geq 1$, consider the map
$s_{ij}^{(k)} : \Gamma_G(k) \rightarrow \mathrm{gr}^k(\mathfrak{J}_G)$ defined by
\[ x \mapsto s_{ij}(x). \]
By Lemma {\rmfamily \ref{L-Pine}}, $s_{ij}^{(k)}$ naturally induces the homomorphism $\mathcal{L}_G(k) \rightarrow \mathrm{gr}^k(\mathfrak{J}_G)$,
which is also denoted by $s_{ij}^{(k)}$ by abuse of language.

\begin{pro}\label{P-blossom}
Assume $\bigcap_{k \geq 1} \Gamma_G(k) = \{1 \}$. Let $k$ be a positive integer. For any $1 \leq p \leq k$, assume
\[ \bigcap_{1 \leq i, j \leq m} \mathrm{Ker}(s_{ij}^{(p)}) = \{ 0 \}. \]
Then $\mathcal{D}_G^m(k) \subset \mathcal{A}_G(k)$.
\end{pro}
\textit{Proof.}
Assume that there exists some $\sigma \in \mathcal{D}_G^m(k)$ such that $\sigma \notin \mathcal{A}_G(k)$.
Since $\mathcal{D}_G^m(k) \subset \mathcal{D}_G^m(1) =\mathcal{A}_G(1)$,
there exists some $1 \leq p \leq k-1$ such that $\sigma \in \mathcal{A}_G(p)$ and
$\sigma \notin \mathcal{A}_G(p+1)$. Thus, we have some $x \in G$ such that
$x^{-1} x^{\sigma} \in \Gamma_G(p+1)$ and $x^{-1} x^{\sigma} \notin \Gamma_G(p+2)$.
By the assumption, for some $1 \leq i, j \leq m$, we see
$s_{ij}(x^{-1} x^{\sigma})$ does not belong to $\mathfrak{J}_G^{p+2}$.
If not, we have $x^{-1} x^{\sigma} \in \bigcap_{1 \leq i, j \leq m} \mathrm{Ker}(s_{ij}^{(p+1)}) = \{ 0 \}$,
and $x^{-1} x^{\sigma} \in \Gamma_G(p+2)$.
Set $\gamma := x^{-1} x^{\sigma} \in \Gamma_G(p+1)$, then
\[\begin{split}
   s_{ij}(x^{\sigma}) & = s_{ij}(x \gamma) \\
   & = s_{ij}(x) + s_{ij}(\gamma) + \sum_{h=1}^m s_{ih}(x) s_{hj}(\gamma) \\
   & \equiv s_{ij}(x) + s_{ij}(\gamma) \pmod{\mathfrak{J}_{F_n}^{p+2}}.
  \end{split}\]
On the other hand, since $\sigma \in \mathcal{D}_G^m(k)$, we have $s_{ij}(x^{\sigma}) - s_{ij}(x) \in \mathfrak{J}_G^{k+1} \subset \mathfrak{J}_G^{p+2}$.
This is the contradiction. Therefore $\sigma \in \mathcal{A}_G(k)$.
$\square$

\vspace{0.5em}

Finally, we introduce generalized Johnson homomorphisms to study the graded quotients $\mathrm{gr}^k(\mathcal{D}_G^m) := \mathcal{D}_G^m(k)/\mathcal{D}_G^m(k+1)$.
To begin with, we consider the action of $\mathrm{Aut}\,G$. 
Since each $\mathcal{D}_G^m(k)$ is a normal subgroup of $\mathrm{Aut}\,G$, the group $\mathrm{Aut}\,G$ naturally acts on 
$\mathrm{gr}^k(\mathcal{D}_G^m)$ by the conjugation from the right. Furthermore since $\mathcal{D}_G^m(1) \supset \mathcal{D}_G^m(2) \supset \cdots$
is a central filtration, the action of $\mathcal{D}_G^m(1)$ on $\mathrm{gr}^k(\mathcal{D}_G^m)$ is trivial. Hence we can consider
each $\mathrm{gr}^k(\mathcal{D}_G^m)$ as an $\mathrm{Aut}\,G/\mathcal{D}_G^m(1)$-module.

\vspace{0.5em}

For any $k \geq 1$ and $\sigma \in \mathcal{D}_G^m(k)$, define the map
$\tilde{\eta}_k(\sigma) : \mathrm{gr}^1(\mathfrak{J}_G) \rightarrow \mathrm{gr}^{k+1}(\mathfrak{J}_G)$ to be
\[ \tilde{\eta}_k(\sigma)(f \pmod{\mathfrak{J}_G^{2}}) := s_{\sigma}(f) \pmod{\mathfrak{J}_G^{k+2}}
    \in \mathrm{gr}^{k+1}(\mathfrak{J}_G) \]
for any $f \in \mathfrak{J}_G$. 
The well-definedness of the map $\tilde{\eta}_k(\sigma)$ follows from Lemma {\rmfamily \ref{L-Peach}}.
It is easily seen that $\tilde{\eta}_k(\sigma)$ is a homomorphism.
Then we have the map $\tilde{\eta}_k : \mathcal{D}_G^m(k) \rightarrow \mathrm{Hom}_{\Q}(\mathrm{gr}^1(\mathfrak{J}_G), \mathrm{gr}^{k+1}(\mathfrak{J}_G))$
defined by $\sigma \mapsto \tilde{\eta}_k(\sigma)$. For any $\sigma, \tau \in \mathcal{D}_G^m(k)$, from (1) of Lemma {\rmfamily \ref{L-Rouge}},
and from Lemma {\rmfamily \ref{L-Peach}}, we see
\[ s_{\sigma \tau}(f) =(s_{\sigma}(f))^{\tau} + s_{\tau}(f) \equiv s_{\sigma}(f) + s_{\tau}(f) \pmod{\mathfrak{J}_G^{k+2}}. \]
This shows that $\tilde{\eta}_k$ is a homomorphism. By the definition, the kernel of $\tilde{\eta}_k$ is $\mathcal{D}_G^m(k+1)$.
Thus $\tilde{\eta}_k$ induces the injective homomorphism
\[ \eta_k : \mathrm{gr}^k(\mathcal{D}_G^m) \rightarrow \mathrm{Hom}_{\Q}(\mathrm{gr}^1(\mathfrak{J}_G), \mathrm{gr}^{k+1}(\mathfrak{J}_G)). \]
The homomorphism $\eta_k$ is an $\mathrm{Aut}\,G/\mathcal{D}_G^m(1)$-equivariant homomorphism.
By using the homomorphisms $\eta_k$, we can consider $\mathrm{gr}^k(\mathcal{D}_G^m)$ as a submodule
of $\mathrm{Hom}_{\Q}(\mathrm{gr}^1(\mathfrak{J}_G), \mathrm{gr}^{k+1}(\mathfrak{J}_G))$, and hence we obtain
\begin{enumerate}
\item Each of $\mathrm{gr}^k(\mathcal{D}_G^m)$ is torsion-free.
\item $\mathrm{dim}_{\Q} (\mathrm{gr}^k(\mathcal{D}_G^m) \otimes_{\Z} \Q) < \infty$.
\end{enumerate}
If $\mathcal{D}_G^m(k) = \mathcal{A}_G(k)$, then the above facts immediately follows from Andreadakis's result for the Andreadakis-Johnson
filtration in \cite{And}.

\section{Extensions of homomorphisms}\label{S-Ext-J}

In this section, we consider to extend the homomorphism $\tilde{\eta}_1$ to $\mathrm{Aut}\,G$ as a crossed homomorphism.
The method we introduce here is based on our previous works in \cite{S19} and \cite{S22}.
Both of them are much inspired by Morita's idea in \cite{Mo5}.

\subsection{Construction of crossed homomorphisms}\label{Ss-cst}
\hspace*{\fill}\ 

\vspace{0.5em}

To begin with, according to the usual convention in homological algebra,
for any group $G$ and $G$-module $M$, we consider that $G$ acts on $M$ from the left if otherwise noted. 
For an associative algebra $\mathcal{R}$, we denote by $\mathrm{Aut}_{\mathrm{(Alg)}}(\mathcal{R})$ the algebra automorphism group of $R$.
For a $\Q$-vector space $M$, we denote by $\mathrm{Aut}\,M$ the group of all $\Q$-linear isomorphisms on $M$.

\vspace{0.5em}

Let $\mathcal{R}$ be a $\Q$-algebra generated by $t_1, \ldots, t_n$, and
$\mathcal{K}$ the ideal of $\mathcal{R}$ generated by $t_1, \ldots, t_n$. We assume that $\mathcal{R} \neq \mathcal{K}$.
Then we have a descending filtration $\mathcal{K} \supset \mathcal{K}^2 \supset \cdots$ of $\mathcal{R}$.
Remark that each $\mathcal{K}/\mathcal{K}^k$ is not $\Q$-algebra, but has a multiplication induced from the quotient $\Q$-algebra
$\mathcal{R}/\mathcal{K}^k$.
For any $f \in \mathcal{K}$, we denote the coset class of $f$ in $\mathcal{K}/\mathcal{K}^k$ by $[f]_k$.
For any $k \geq 2$, set
\[ \overline{\mathrm{Aut}}(\mathcal{K}/\mathcal{K}^k) := \{ \sigma \in \mathrm{Aut}(\mathcal{K}/\mathcal{K}^k) \,|\, \sigma([\gamma \gamma']_k)
      = \sigma([\gamma]_k) \, \sigma([\gamma']_k), \,\,\,
         \gamma, \gamma' \in \mathcal{K} \}. \]
Since
$\mathcal{K}/\mathcal{K}^k \cong \mathcal{K}/\mathcal{K}^2 \oplus \mathcal{K}^2/\mathcal{K}^3 \oplus \cdots \oplus \mathcal{K}^{k-1}/\mathcal{K}^k$
as a $\Q$-vector space,
we see $\mathrm{Aut}(\mathcal{K}/\mathcal{K}^k) = \mathrm{GL}_{\Q}(\mathcal{K}/\mathcal{K}^k) \cong \mathrm{GL}(m,\Q)$ where
$m=\sum_{i=1}^{k-1} \mathrm{dim}_{\Q}(\mathcal{K}^i/\mathcal{K}^{i+1})$. 
Note that $\overline{\mathrm{Aut}}(\mathcal{K}/\mathcal{K}^2)=\mathrm{Aut}(\mathcal{K}/\mathcal{K}^2)=\mathrm{GL}_{\Q}(\mathcal{K}/\mathcal{K}^2)$.
In general, $\overline{\mathrm{Aut}}(\mathcal{K}/\mathcal{K}^k)$
is not equal to $\mathrm{Aut}(\mathcal{K}/\mathcal{K}^k)$. 

\vspace{0.5em}

First, we consider an embedding of $\mathrm{Hom}_{\Q}(\mathrm{gr}^1(\mathcal{K}), \mathrm{gr}^2(\mathcal{K}))$
into $\overline{\mathrm{Aut}}(\mathcal{K}/\mathcal{K}^3)$.
For any $f \in \mathrm{Hom}_{\Q}(\mathrm{gr}^1(\mathcal{K}), \mathrm{gr}^2(\mathcal{K}))$, define the map
$\widetilde{f} : \mathcal{K}/\mathcal{K}^3 \rightarrow \mathcal{K}/\mathcal{K}^3$
by
\[ \widetilde{f}([\gamma]_3) := [\gamma]_3 + f([\gamma]_2) \]
for any $\gamma \in \mathcal{K}$.

\begin{pro}\label{L-cofley}
With the above notation, for any $f \in \mathrm{Hom}_{\Q}(\mathrm{gr}^1(\mathcal{K}), \mathrm{gr}^2(\mathcal{K}))$,
we see $\widetilde{f} \in \overline{\mathrm{Aut}}(\mathcal{K}/\mathcal{K}^3)$, and the map
\[ \iota : \mathrm{Hom}_{\Q}(\mathrm{gr}^1(\mathcal{K}), \mathrm{gr}^2(\mathcal{K})) \rightarrow \overline{\mathrm{Aut}}\,(\mathcal{K}/\mathcal{K}^3) \]
defined by $f \mapsto \widetilde{f}$ is injective.
\end{pro}
\textit{Proof.}
First, we show $\widetilde{f}$ is a homomorphism. For any $\gamma, \gamma' \in \mathcal{K}$,
\[\begin{split}
   \widetilde{f}([\gamma]_3 + [\gamma']_3) & = \widetilde{f}([\gamma + \gamma']_3)
         = [\gamma + \gamma']_3 + f([\gamma + \gamma']_2) \\
       & = ([\gamma]_3 + f([\gamma]_2) ) + ([\gamma']_3 + f([\gamma']_2) ) \\
       & = \widetilde{f}([\gamma]_3) + \widetilde{f}([\gamma']_3).
  \end{split}\]
Thus $\widetilde{f}$ is a homomorphism. Furthermore, $\widetilde{f}$ satisfies
\[\begin{split}
   \widetilde{f}([\gamma]_3 [\gamma']_3) & = \widetilde{f}([\gamma \gamma']_3) = [\gamma \gamma']_3  + f([\gamma \gamma']_2) = [\gamma]_3 [\gamma']_3 \\
    & = ([\gamma]_3 + f([\gamma]_2))([\gamma']_3 + f([\gamma']_2)) = \widetilde{f}([\gamma]_3) \widetilde{f}([\gamma']_3)
  \end{split}\]
for any $\gamma, \gamma' \in \mathcal{K}$.
On the other hand, we have $\widetilde{f+g}=\widetilde{f} \circ \tilde{g}$ for any $f, g \in \mathrm{Hom}_{\Q}(\mathrm{gr}^1(\mathcal{K}), \mathrm{gr}^2(\mathcal{K}))$.
In fact, for any $\gamma \in \mathcal{K}$,
\[\begin{split}
   (\widetilde{f+g})([\gamma]_3) & = [\gamma]_3 + f([\gamma]_2) +g([\gamma]_2), \\
   (\widetilde{f} \circ \tilde{g})([\gamma]_3) & = \widetilde{f}([\gamma]_3 + g([\gamma]_2))
     = [\gamma]_3 + g([\gamma]_2) + f([\gamma]_2).
  \end{split}\]
This shows that $\iota$ is a homomorphism.
On the other hand, for the zero map $0 \in \mathrm{Hom}_{\Q}(\mathrm{gr}^1(\mathcal{K}), \mathrm{gr}^2(\mathcal{K}))$, it is obvious that
$\widetilde{0} = \mathrm{id}_{\mathcal{K}/\mathcal{K}^3}$. Hence each $\widetilde{f}$ has its inverse map $\widetilde{f}^{-1} = \widetilde{-f}$.
This means $\widetilde{f}$ an automorphism on $\mathcal{K}/\mathcal{K}^3$.

\vspace{0.5em}

Finally, we show that $\iota$ is injective. Assume that $\widetilde{f}=\mathrm{id}_{\mathcal{K}/\mathcal{K}^3}$ for some
$f \in \mathrm{Hom}_{\Q}(\mathrm{gr}^1(\mathcal{K})$, $\mathrm{gr}^2(\mathcal{K}))$. Then for any $[\gamma]_3 \in \mathcal{K}/\mathcal{K}^3$, we have
\[ \widetilde{f}([\gamma]_3) = [\gamma]_3 + f([\gamma]_2) = [\gamma]_3. \]
Hence we obtain $f([\gamma]_2)=0$ for any $[\gamma]_2 \in \mathcal{K}/\mathcal{K}^2$, and $f=0$. This shows $\iota$ is injective. $\square$

\vspace{0.5em}

Let
\[ \varphi : \overline{\mathrm{Aut}}(\mathcal{K}/\mathcal{K}^3) \rightarrow \overline{\mathrm{Aut}}(\mathcal{K}/\mathcal{K}^2) \]
be the natural homomorphism induced from the projection $\mathcal{K}/\mathcal{K}^3 \rightarrow \mathcal{K}/\mathcal{K}^2$.
\begin{pro}\label{P-spt}
The sequence 
\begin{equation}
 0 \rightarrow \mathrm{Hom}_{\Q}(\mathrm{gr}^1(\mathcal{K}), \mathrm{gr}^2(\mathcal{K})) \xrightarrow{\iota}
   \overline{\mathrm{Aut}}\,(\mathcal{K}/\mathcal{K}^3) \xrightarrow{\varphi} \mathrm{Im}(\varphi) \rightarrow 1 \label{eq-exact}
\end{equation}
is a split group extension.
\end{pro}
\textit{Proof.}
First, we show the above sequence is exact. Namely, it suffices to show
$\mathrm{Im}(\iota)=\mathrm{Ker}(\varphi)$. The fact that $\mathrm{Im}(\iota) \subset \mathrm{Ker}(\varphi)$
follows from
\[\begin{split}
   (\varphi \circ \iota)(f)([\gamma]_2) & = \varphi(\widetilde{f})([\gamma]_2) = [\gamma]_3 + f([\gamma]_2) \pmod{\mathcal{K}^2} \\
     & = [\gamma]_2
  \end{split}\]
for any $f \in \mathrm{Hom}_{\Q}(\mathrm{gr}^1(\mathcal{K}), \mathrm{gr}^2(\mathcal{K}))$ and $\gamma \in \mathcal{K}$.
To show $\mathrm{Im}(\iota) \supset \mathrm{Ker}(\varphi)$, take any $f \in \mathrm{Ker}(\varphi)$.
Define $g \in \mathrm{Hom}_{\Q}(\mathrm{gr}^1(\mathcal{K}), \mathrm{gr}^2(\mathcal{K}))$ to be
\begin{equation}\label{eq-aoi}
 g([\gamma]_2) := f([\gamma]_3) -[\gamma]_3
\end{equation}
for any $\gamma \in \mathcal{K}$. The map $g$ is well-defined.
In fact, for any $\gamma, \gamma' \in \mathcal{K}$ such that $[\gamma]_2=[\gamma']_2$, if we set
$\gamma' - \gamma = \varepsilon \in \mathcal{K}^2$, then we have
\[\begin{split}
   g([\gamma']_2) &= f([\gamma']_3) -[\gamma']_3 = f([\gamma + \varepsilon]_3) -[\gamma + \varepsilon]_3 \\
                  &= f([\gamma]_3) -[\gamma]_3 + f([\varepsilon]_3) -[\varepsilon]_3 = f([\gamma]_3) -[\gamma]_3 \\
                  &= g([\gamma]_2).
  \end{split}\]
Here we remark that $f([\varepsilon]_3) -[\varepsilon]_3 =0 \in \mathcal{K}^2/\mathcal{K}^3$ since $\varepsilon \in \mathcal{K}^2$ and $f \in \mathrm{Ker}(\varphi)$.
It is easy to show that $g$ is a homomorphism. Furthermore, for any $\gamma \in \mathcal{K}$,
\[ \widetilde{g}([\gamma]_3) = [\gamma]_3 + g([\gamma]_2) = f([\gamma]_3). \]
This shows $f=\widetilde{g}=\iota(g) \in \mathrm{Im}(\iota)$.

\vspace{0.5em}

Next, we construct the section $s : \mathrm{Im}(\varphi) \rightarrow \overline{\mathrm{Aut}}\,(\mathcal{K}/\mathcal{K}^3)$ of (\ref{eq-exact}).
Since $t_1, \ldots, t_n \pmod{\mathcal{K}^2}$ generate $\mathrm{gr}^1(\mathcal{K})$ as a $\Q$-vector space,
we can chose $t_{i_1}, \ldots t_{i_p}$ such that $t_{i_1}, \ldots t_{i_p} \pmod{\mathcal{K}^2}$ is a basis of $\mathrm{gr}^1(\mathcal{K})$.
Fix it, and write $\gamma_j := t_{i_j}$ for simplicity. Furthermore, 
elements $t_{i_j} t_{i_l} \pmod{\mathcal{K}^3}$ for $1 \leq j, l \leq p$ generate $\mathrm{gr}^2(\mathcal{K})$ as a $\Q$-vector space,
we can chose elements from them such that their coset classes form a basis of $\mathrm{gr}^2(\mathcal{K})$.
Say them $\gamma_{p+1}, \ldots, \gamma_{p+q} \in \mathcal{K}^2$. Then
$[\gamma_1]_3, \ldots [\gamma_p]_3, [\gamma_{p+1}]_3, \ldots, [\gamma_{p+q}]_3$ is a basis of $\mathcal{K}/\mathcal{K}^3$ as a $\Q$-vector space.

\vspace{0.5em}

For any $\beta \in \mathrm{Im}(\varphi)$, there exists an element
$\beta' \in \overline{\mathrm{Aut}}\,(\mathcal{K}/\mathcal{K}^3)$ such that
$\varphi(\beta') = \beta$. Since $\beta'$ preserve the multiplication, $\beta'$ is uniquely determined by the images of $[\gamma_j]_3$
for all $1 \leq j \leq p$. For any $1 \leq j \leq p$, set
\[ \beta'([\gamma_j]_3) := c_{1j} [\gamma_1]_3 + \cdots + c_{pj}[\gamma_p]_3 
     + c_{p+1, j} [\gamma_{p+1}]_3 + \cdots + c_{p+q, j} [\gamma_{p+q}]_3 \]
for some $c_{ij} \in \Q$. Since $\beta \in \overline{\mathrm{Aut}}(\mathcal{K}/\mathcal{K}^2)=\mathrm{GL}_{\Q}(\mathcal{K}/\mathcal{K}^2)$,
if we set
\[ v_j := c_{1j} [\gamma_1]_2 + \cdots + c_{pj}[\gamma_p]_2 \]
for any $1 \leq j \leq p$, then $v_1, v_2, \ldots, v_p$ is a $\Q$-linear basis of $\mathrm{gr}^1(\mathcal{K})=\mathcal{K}/\mathcal{K}^2$.
Let $\delta = \delta_{\beta'} : \mathrm{gr}^1(\mathcal{K}) \rightarrow \mathrm{gr}^2(\mathcal{K})$ be the $\Q$-linear map given by
\[ \delta(v_j) = - (c_{p+1, j} [\gamma_{p+1}]_3 + \cdots + c_{p+q, j} [\gamma_{p+q}]_3) \]
for any $1 \leq j \leq p$. Then we obtain
\[\begin{split}
  (\iota(\delta) \circ \beta')([\gamma_j]_3)
   & = \iota(\delta)(c_{1j} [\gamma_1]_3 + \cdots + c_{pj}[\gamma_p]_3 + c_{p+1, j} [\gamma_{p+1}]_3 + \cdots + c_{p+q, j} [\gamma_{p+q}]_3) \\
   & = c_{1j} [\gamma_1]_3 + \cdots + c_{pj}[\gamma_p]_3 + c_{p+1, j} [\gamma_{p+1}]_3 + \cdots + c_{p+q, j} [\gamma_{p+q}]_3 \\
   & \hspace{2em} + \delta(c_{1j} [\gamma_1]_2 + \cdots + c_{pj}[\gamma_p]_2) \\
   & = c_{1j} [\gamma_1]_3 + \cdots + c_{pj}[\gamma_p]_3 + c_{p+1, j} [\gamma_{p+1}]_3 + \cdots + c_{p+q, j} [\gamma_{p+q}]_3 \\
   & \hspace{2em} + \delta(v_j) \\
   & = c_{1j} [\gamma_1]_3 + \cdots + c_{pj}[\gamma_p]_3
  \end{split}\]
for any $1 \leq j \leq p$.
Consider the map $s : \overline{\mathrm{Aut}}\,(\mathcal{K}/\mathcal{K}^2) \rightarrow \overline{\mathrm{Aut}}\,(\mathcal{K}/\mathcal{K}^3)$
define by $\beta \mapsto \iota(\delta) \circ \beta'$.
We can see that $s$ is a homomorphism and is the required section. Hence the exact sequence (\ref{eq-exact}) splits.
$\square$

\vspace{0.5em}

Observe the following lemma.

\begin{lem}\label{L-ext}
Let
\[ 0 \rightarrow K \rightarrow G \rightarrow N \rightarrow 1 \]
be a split extension of groups over $N$ with $K$ additive abelian group. For any $g$, there exist unique elements $k_g \in K$ and $n_g \in N$
such that $g=k_g n_g$. Then the map $k : G \rightarrow K$ defined by $g \mapsto k_g$ is a crossed homomorphism.
\end{lem}
\noindent
Since the proof is an easy exercise, we leave it to the leader.
By applying this lemma to the split extension (\ref{eq-exact}),
we obtain the crossed homomorphism
\begin{equation}\label{eq-sushi}
 \theta : \overline{\mathrm{Aut}}(\mathcal{K}/\mathcal{K}^3) \rightarrow \mathrm{Hom}_{\Q}(\mathrm{gr}^1(\mathcal{K}), \mathrm{gr}^2(\mathcal{K})).
\end{equation}

\subsection{Extensions of $\tilde{\eta}_1$}\label{Ss-Ext}
\hspace*{\fill}\ 

\vspace{0.5em}

Here, we show that $\tilde{\eta}_1$ can be extended to $\mathrm{Aut}\,G$ as a crossed homomorphism.
According to the usual convention in homological algebra,
for any group $G$ and $G$-module $M$, we consider that $G$ acts on $M$ from the left if otherwise noted. Hence, the right actions
mentioned above are read as the left one in the natural way. For example, for any $\sigma \in \mathrm{Aut}\,G$ and $x \in G$,
the left action of $\sigma$ on $s_{ij}(x)$ is given by
\[ \sigma \cdot s_{ij}(x) = s_{ij}(x^{\sigma^{-1}}). \]

\vspace{0.5em}

Consider the associative algebra $\mathfrak{R}_{\Q}^m(G)$ and the filtration
$\mathfrak{J}_G \supset \mathfrak{J}_G^2 \supset \mathfrak{J}_G^3 \supset \cdots$.
By using (\ref{eq-sushi}), we obtain the crossed homomorphism
\[ \theta : \overline{\mathrm{Aut}}(\mathfrak{J}_G/\mathfrak{J}_G^3)
    \rightarrow \mathrm{Hom}_{\Q}(\mathrm{gr}^1(\mathfrak{J}_G), \mathrm{gr}^2(\mathfrak{J}_G)). \]
Let $\mathrm{Aut}\,G \rightarrow \mathrm{Aut}(\mathfrak{J}_G/\mathfrak{J}_G^3)$ be the natural homomorphism induced from the action of
$\mathrm{Aut}\,G$ on $\mathfrak{J}_G/\mathfrak{J}_G^3$. Then its image is contained in $\overline{\mathrm{Aut}}(\mathfrak{J}_G/\mathfrak{J}_G^3)$.
By composing $\theta$ and this natural homomorphism,
we obtain the crossed homomorphism
\[ \theta_{G} : \mathrm{Aut}\,G \rightarrow \mathrm{Hom}_{\Q}(\mathrm{gr}^1(\mathfrak{J}_G), \mathrm{gr}^2(\mathfrak{J}_G)). \]
By observing (\ref{eq-aoi}), we can see that $\theta_{\mathcal{D}}$ is an extension of
$\tilde{\eta}_1 : \mathcal{D}_G^m(1) \rightarrow \mathrm{Hom}_{\Q}(\mathrm{gr}^1(\mathfrak{J}_G), \mathrm{gr}^{2}(\mathfrak{J}_G))$.

\section{On the $\mathrm{SL}(m,\C)$-representation algebra of $F_n$}\label{S-Free}

\vspace{0.5em}

In this section, we study the $\mathrm{SL}(m,\C)$-representation algebra for free groups.
Let $F_n$ be the free group on rank $n$ with basis $x_1, x_2, \ldots, x_n$.

\subsection{The graded quotients $\mathrm{gr}^k(\mathfrak{J}_{F_n})$}\label{Ss-FRep}
\hspace*{\fill}\ 

\vspace{0.5em}

First, we give a basis of $\mathrm{gr}^k(\mathfrak{J}_{F_n})$ as a $\Q$-vector space.
For any $k \geq 1$, recall that (the coset classes of) the elements in
\[ T_k := \Big{\{} \prod_{\substack{ 1 \leq i, j \leq m \\[1pt] (i,j) \neq (m,m)}} \prod_{l=1}^n s_{ij}(x_l)^{e_{ij,l}}
    \,\Big{|}\, e_{ij,l} \geq 0, \,\,\,
    \sum_{\substack{ 1 \leq i, j \leq m \\[1pt] (i,j) \neq (m,m)}} \sum_{l=1}^n e_{ij,l} = k \Big{\}}
   \subset \mathfrak{J}_{F_n}^k \] 
generate $\mathrm{gr}^k(\mathfrak{J}_{F_n})$ as a $\Q$-vector space.
\begin{pro}\label{P-Rythm}
For each $k \geq 1$, the set $T_k \pmod{\mathfrak{J}_{F_n}^{k+1}}$ forms a basis of $\mathrm{gr}^k(\mathfrak{J}_{F_n})$ as a $\Q$-vector space.
\end{pro}
\textit{Proof.}
In order to show that $T_k \pmod{\mathfrak{J}_{F_n}^{k+1}}$ is linearly independent over $\Q$, assume
\[ {\sum}' a(e_{ij,l}) 
    \prod_{\substack{ 1 \leq i, j \leq m \\[1pt] (i,j) \neq (m,m)}} \prod_{l=1}^n s_{ij}(x_l)^{e_{ij,l}}
    \equiv 0 \pmod{\mathfrak{J}_{F_n}^{k+1}} \]
where the above sum $\sum'$ runs over all tuples $(e_{ij,l})$ for $1 \leq i, j \leq m$, $(i,j) \neq (m,m)$ and $1 \leq l \leq n$
such that the sum of $e_{ij,l}$ is equal to $k$.
Denote by $f \in \mathfrak{R}_{\Q}^m(F_n)$ the left hand side of the above equation, and assume $f \in \mathfrak{J}_{F_n}^{k+1}$.
In order to obtain $a(e_{ij,l})=0$, we construct some $\mathrm{SL}(m,\C)$-representations of $F_n$, and investigate relations
among $a(e_{ij,l})$ by observing the values of them by $f : R(F_n) \rightarrow \C$.

\vspace{0.5em}

For any $n \geq 1$, consider
\[\begin{split}
   D_m & := \{ (z_{ij})_{1 \leq i, j \leq m} \in \C^{m^2} \,|\, \det{((z_{ij})_{1 \leq i, j \leq m})} \neq 0 \} \subset \C^{m^2}.
  \end{split}\]
The set $D_m$ is an open subset in $\C^{m^2}$.
For any $1 \leq l \leq n$, take any $(z_{ij,l})_{1 \leq i, j \leq m-1} \in D_{m-1}$, and 
any $z_{im,l}, z_{mj, l} \in \C$ for $1 \leq i, j \leq m-1$.
Define the representation $\rho : F_n \rightarrow \mathrm{SL}(m,\C)$ by
\[ \rho(x_l) := (z_{ij,l})_{1 \leq i, j \leq m} \]
for any $1 \leq l \leq n$. Here, $z_{mm, l}$ is defined by
\[\begin{split}
   z_{mm, l} := Z_{mm,l}^{-1}(1 - (z_{1m,l} Z_{1m,l} + \cdots + z_{m-1\, m, l} Z_{m-1\, m, l})), 
  \end{split}\]
and $Z_{ij,l}$ is the $(i,j)$-cofactor of the matrix $(z_{ij,l})_{1 \leq i, j \leq m}$. Note that
\[ Z_{mm,l} = \det{((z_{ij,l})_{1 \leq i, j \leq m-1})} \neq 0. \]

\vspace{0.5em}

Then we have
\[ f(\rho) = {\sum}' a(e_{ij,l}) 
    \prod_{\substack{ 1 \leq i, j \leq m \\[1pt] (i,j) \neq (m,m)}} \prod_{l=1}^n (\overline{z}_{ij,l})^{e_{ij,l}} \]
where
\[ \overline{z}_{ij,l} = \begin{cases}
                          z_{ij,l}, \hspace{1em} & \text{if} \hspace{0.5em} i \neq j, \\
                          z_{ii,l} -1, & \text{if} \hspace{0.5em} i = j.
                         \end{cases}\]
On the other hand, since $f \in \mathfrak{J}_{F_n}^{k+1}$, we see that
$f(\rho)$ can be written as a polynomial among $\overline{z}_{ij,l}$ with degree $\geq k+1$.
Since we can take $(z_{ij,l})_{1 \leq i, j \leq m-1} \in D_{m-1}$, and $z_{im,l}, z_{mj, l} \in \C$ arbitrary, 
by the uniqueness of the expression, we obtain
$a(e_{ij,l})=0$ for any tuple $(e_{ij,l})$.
$\square$

\vspace{0.5em}

As a corollary, we see
\[ \mathrm{dim}_{\Q} \mathrm{gr}^k(\mathfrak{J}_{F_n}) = \binom{(m^2-1)n+k-1}{k}. \]

\vspace{0.5em}

Here we give a few problems.
\begin{prob}
\begin{enumerate}
\item[(1)] Is $\bigcap_{k \geq 1} \mathfrak{J}_{F_n}^k$ is equal to $\{ 0 \}$?
\item[(2)] Is the algebra $\mathfrak{R}_{\Q}(F_n)$ isomorphic to $\mathrm{SL}_m$-universal representation algebra of $F_n$ over $\Q$?
\end{enumerate}
\end{prob}
\noindent
In \cite{S20}, we gave affirmative answers to the above questions for $m=2$.
It seems, however, quite difficult to apply the same method to the $\mathrm{SL}_m$-representation case for $m \geq 3$ due to the combinatorial complexity.

\vspace{0.5em}

Now, we review $\mathrm{Aut}\,F_n$ and the IA-automorphism group $\mathrm{IA}(F_n)=\mathcal{A}_{F_n}(1)$ of $F_n$.
Let $P$, $Q$, $S$ and $U$ be the automorphisms of $F_n$ given by specifying its images of the basis $x_1, \ldots, x_n$ as follows:

\vspace{0.5em}

\begin{center}
\begin{tabular}{|c|c|c|c|c|c|c|} \hline
           & $x_1$      & $x_2$ & $x_3$ & $\cdots$ & $x_{n-1}$ & $x_n$ \\ \hline
  $P$      & $x_2$      & $x_1$ & $x_3$ & $\cdots$ & $x_{n-1}$ & $x_n$ \\ 
  $Q$      & $x_2$      & $x_3$ & $x_4$ & $\cdots$ & $x_{n}$   & $x_1$ \\ 
  $S$      & $x_1^{-1}$ & $x_2$ & $x_3$ & $\cdots$ & $x_{n-1}$ & $x_n$ \\ 
  $U$      & $x_1 x_2$  & $x_2$ & $x_3$ & $\cdots$ & $x_{n-1}$ & $x_n$ \\ \hline
\end{tabular}
\end{center}

\vspace{0.5em}

\noindent
In 1924, Nielsen \cite{Ni1} showed that $\mathrm{Aut}\,F_n$ is generated by $P$, $Q$, $S$ and $U$, and gave finitely many relators among them.
On the other hand, Magnus \cite{Mag} showed that for any $n \geq 3$, $\mathrm{IA}(F_n)$ is finitely generated by automorphisms
\[ K_{ij} : x_t \mapsto \begin{cases}
               {x_j}^{-1} x_i x_j, & t=i, \\
               x_t,                & t \neq i
              \end{cases}\]
for distinct $1 \leq i, \, j \leq n$, and
\[  K_{ijl} : x_t \mapsto \begin{cases}
               x_i [x_j, x_l], & t=i, \\
               x_t,            & t \neq i
              \end{cases}\] 
for distinct $1 \leq i, \, j, \, l \leq n$ and $j<l$. Let $H$ be the abelianization $F_n/\Gamma_{F_n}(2)$ of $F_n$.
By observing the images of Nielsen's generators, we can see that the natural homomorphism $\mathrm{Aut}\,F_n \rightarrow \mathrm{Aut}\,H$ is surjective.
By fixing the basis of the free abelian group $H$ induced from that $x_1, \ldots, x_n$ of $F_n$, we can identify
$\mathrm{Aut}(H)$ with the general linear group $\mathrm{GL}(n,\Z)$. Hence we have
$\mathrm{Aut}\,F_n/\mathrm{IA}(F_n) \cong \mathrm{GL}(n,\Z)$.

\vspace{0.5em}

In \cite{S20}, we showed that $\mathcal{D}_{F_n}^2(k)=\mathcal{A}_{F_n}(k)$ for $1 \leq k \leq 4$.
Hence, from Theorem {\rmfamily \ref{T-Flower}}, we have the following.
\begin{thm}\label{T-Flower-2}
For any $m \geq 2$, we have $\mathcal{D}_{F_n}^m(k)=\mathcal{A}_{F_n}(k)$ for any $1 \leq k \leq 4$.
\end{thm}
\noindent
In particular, we see that $\mathrm{GL}(n,\Z) \cong \mathrm{Aut}\,F_n/\mathcal{D}_{F_n}^m(1)$ naturally acts on $\mathrm{gr}^k(\mathfrak{J}_{F_n})$
for any $k \geq 1$. Set $H_{\Q} := H \otimes_{\Z} \Q$.
\begin{pro}\label{P-dim-2}
For any $n \geq 2$ and $k \geq 1$, we have
\[ \mathrm{gr}^k(\mathfrak{J}_{F_n}) \cong {\bigoplus}^{'}
    \bigotimes_{\substack{ 1 \leq i, j \leq m \\[1pt] (i,j) \neq (m,m)}} S^{e_{ij}} H_{\Q} \]
as a $\mathrm{GL}(n,\Z)$-module. Here the sum runs over all tuples $(e_{ij})$ for $1 \leq i, j \leq m$ and $(i,j) \neq (m,m)$ such that
the sum of the $e_{ij}$ is equal to $k$.
\end{pro}
\textit{Proof.}
Let $\mathfrak{M}$ be the right hand side of the above equation. First, for any $1 \leq i, j \leq m$ and $e \geq 1$, the homomorphism
$f_{ij}^e : S^e H_{\Q} \rightarrow \mathrm{gr}^e(\mathfrak{J}_{F_n})$ defined by
\[ x_1^{l_1} x_2^{l_2} \cdots x_n^{l_n} \mapsto s_{ij}(x_1)^{l_1} s_{ij}(x_2)^{l_2} \cdots s_{ij}(x_n)^{l_n} \pmod{\mathfrak{J}_{F_n}^{e+1}} \]
for $l_1 + l_2 + \cdots + l_n=e$
is $\mathrm{Aut}\,F_n$-equivariant. In fact, for any Nielsen generators $\sigma = P, Q, S$ and $U$ of $\mathrm{Aut}\,F_n$, we can check
$f_{ij}^e(x^{\sigma}) =(f_{ij}^e(x))^{\sigma}$ for any $x=x_1^{l_1} x_2^{l_2} \cdots x_n^{l_n}$. For example, we see
\[\begin{split}
  f_{ij}^e((x_1^{l_1} x_2^{l_2} \cdots x_n^{l_n})^U) & = f_{ij}^e((x_1+x_2)^{l_1} x_2^{l_2} \cdots x_n^{l_n}) \\
  & = f_{ij}^e \Big{(} \sum_{t = 0}^{l_1} \binom{l_1}{t} x_1^t x_2^{l_1-t}  x_2^{l_2} \cdots x_n^{l_n} \Big{)} \\
  & \equiv \sum_{t = 0}^{l_1} \binom{l_1}{t}  s_{ij}(x_1)^{t} s_{ij}(x_2)^{l_1 -t + l_2} \cdots s_{ij}(x_n)^{l_n} \pmod{\mathfrak{J}_{F_n}^{e+1}} \\
  & \equiv s_{ij}(x_1 x_2)^{l_1} s_{ij}(x_2)^{l_2} \cdots s_{ij}(x_n)^{l_n} \pmod{\mathfrak{J}_{F_n}^{e+1}} \\
  & = (f_{ij}^e(x_1^{l_1} x_2^{l_2} \cdots x_n^{l_n}))^U.
  \end{split}\]

Since $\mathrm{IA}(F_n)$ trivially acts on both of $S^e H_{\Q}$ and $\mathrm{gr}^e(\mathfrak{J}_{F_n})$,
we obtain the surjective $\mathrm{GL}(n,\Z)$-equivariant homomorphism
$F : \mathfrak{M} \rightarrow \mathrm{gr}^k(\mathfrak{J}_{F_n})$ defined by
\[\begin{split}
   {\sum}' & a_{e_{11}, \ldots, e_{m m-1}} \bigotimes_{\substack{ 1 \leq i, j \leq m \\[1pt] (i,j) \neq (m,m)}} X_{e_{ij}} \\
    & \mapsto {\sum}' a_{e_{11}, \ldots, e_{m m-1}} \bigotimes_{\substack{ 1 \leq i, j \leq m \\[1pt] (i,j) \neq (m,m)}} f_{ij}^{e_{ij}}(X_{e_{ij}})
  \end{split}\]
for any $X_{e_{ij}} \in S^{e_{ij}} H_{\Q}$ and $a_{e_{11}, \ldots, e_{m m-1}} \in \Q$.
Here the sum runs over all tuples $(e_{ij})$ for $1 \leq i, j \leq m$ and $(i,j) \neq (m,m)$ such that
the sum of $e_{ij}$ is equal to $k$.
The surjectivity of $F$ follows from Proposition {\rmfamily \ref{P-Rythm}}. In fact, for any element
\[ Y := \prod_{\substack{ 1 \leq i, j \leq m \\[1pt] (i,j) \neq (m,m)}} \prod_{l=1}^n s_{ij}(x_l)^{e_{ij,l}} \] 
in the basis $T_k$ of $\mathrm{gr}^k(\mathfrak{J}_{F_n})$, set
\[ X := {\sum}'' \bigotimes_{\substack{ 1 \leq i, j \leq m \\[1pt] (i,j) \neq (m,m)}}
      x_1^{e_{ij,1}} \cdots x_n^{e_{ij,n}}
   \in \mathfrak{M} \]
where the sum runs over all tuples $(e_{ij,l})$ for $1 \leq i, j \leq m$, $(i,j) \neq (m,m)$ and $1 \leq l \leq n$
such that the sum of $e_{ij, l}$ is equal to $k$.
Then we have $Y=F(X)$.

\vspace{0.5em}

Next, we prove that $F$ is an isomorphism by showing that the dimensions of $\mathfrak{M}$ and $\mathrm{gr}^k(\mathfrak{J}_{F_n})$ as $\Q$-vector spaces are equal.
The basis $T_k$ can be expressed as
\[\begin{split}
   T_k & = \Big{\{} \prod_{\substack{ 1 \leq i, j \leq m \\[1pt] (i,j) \neq (m,m)}} \prod_{l=1}^n s_{ij}(x_l)^{e_{ij,l}}
    \,\Big{|}\, e_{ij,l} \geq 0, \,\,\,
    \sum_{\substack{ 1 \leq i, j \leq m \\[1pt] (i,j) \neq (m,m)}} \sum_{l=1}^n e_{ij,l} = k \Big{\}} \\
      & = \Big{\{} \prod_{\substack{1 \leq i, j \leq m \\[1pt] (i,j) \neq (m,m)}}
       s_{ij}(x_1)^{e_{ij,1}} \cdots s_{ij}(x_n)^{e_{ij,n}} \,\Big{|}\, \sum_{l=1}^n e_{ij,l} =e_{ij}, \hspace{0.5em}
         \sum_{\substack{ 1 \leq i, j \leq m \\[1pt] (i,j) \neq (m,m)}} e_{ij} =k \Big{\}}
  \end{split}\]
From the last term of the above equation, we see that the number of elements in $T_k$ is equal to $\mathrm{dim}_{\Q} \mathfrak{M}$.
This completes the proof of Proposition {\rmfamily \ref{P-dim-2}}.
$\square$

\subsection{On the extension of $\tilde{\eta}_1$}\label{Ss-1-coh}
\hspace*{\fill}\ 

\vspace{0.5em}

Let us consider the homomorphism $\tilde{\eta}_1$ and its extension to $\mathrm{Aut}\,F_n$ as a crossed homomorphism.
Observe
\[\begin{split}
   s_{ij}(x[y,z]) - s_{ij}(x) & = s_{ij}(x) + s_{ij}([y,z]) + \sum_{k=1}^m s_{ik}(x)s_{kj}([y,z]) - s_{ij}(x) \\
    & \equiv s_{ij}([y,z]) \pmod{\mathfrak{J}_{F_n}^3} \\
    & \equiv \sum_{k=1}^m (s_{ik}(y) s_{kj}(z) - s_{ik}(z) s_{kj}(y)) \pmod{\mathfrak{J}_{F_n}^3}
  \end{split}\]
for any $x, y, z \in F_n$. The last equality is obtained from the straightforward calculation with (\ref{eq-pooh}).
By using this,
we can easily calculate the images of Magnus generators $K_{ij}$ and $K_{ijl}$ of $\mathrm{IA}(F_n)$ by $\tilde{\eta}_1$ as follows.
From Proposition {\rmfamily \ref{P-Rythm}}, we see that elements $s_{ij}(x_l)$ for $(i,j) \neq (m,m)$ and $1 \leq l \leq n$ form a basis of $\mathrm{gr}^1(\mathfrak{J}_{F_n})$.
Let $s_{ij}(x_l)^*$ be its dual basis of $\mathrm{gr}^1(\mathfrak{J}_{F_n})^* := \mathrm{Hom}_{\Q}(\mathrm{gr}^1(\mathfrak{J}_{F_n}), \Q)$.
Then we have
\begin{equation}
 \begin{split}
   \tilde{\eta}_1(K_{ij})
    & = \sum_{\substack{1 \leq p, q \leq m \\[1pt] (p,q) \neq (m,m)}} \sum_{k=1}^m
       s_{pq}(x_i)^* \otimes (s_{pk}(x_i) s_{kq}(x_j) - s_{pk}(x_j) s_{kq}(x_i)), \\
   \tilde{\eta}_1(K_{ijl}) & = \sum_{\substack{1 \leq p, q \leq m \\[1pt] (p,q) \neq (m,m)}} \sum_{k=1}^m
       s_{pq}(x_i)^* \otimes (s_{pk}(x_j) s_{kq}(x_l) - s_{pk}(x_l) s_{kq}(x_j)).
 \end{split} \label{eq-image-1}
\end{equation}

From the argument in Subsection {\rmfamily \ref{Ss-Ext}}, we have the crossed homomorphism
\[ \theta_{F_n} : \mathrm{Aut}\,F_n \rightarrow \mathrm{Hom}_{\Q}(\mathrm{gr}^1(\mathfrak{J}_{F_n}), \mathrm{gr}^2(\mathfrak{J}_{F_n})), \]
which is an extension of $\tilde{\eta}_1$ to $\mathrm{Aut}\,F_n$.
In this subsection, we study this crossed homomorphism, and relations to $f_K$ and $f_M$.

\vspace{0.5em}

To begin with, we determine the images of Nielsen's generators by $\theta_{F_n}$.
Let $\rho_k : \mathrm{Aut}\,F_n \rightarrow \overline{\mathrm{Aut}}(\mathfrak{J}_{F_n}/\mathfrak{J}_{F_n}^k)$
be the homomorphism induced from the action of $\mathrm{Aut}\,F_n$ on $\mathfrak{J}_{F_n}/\mathfrak{J}_{F_n}^k$.
By direct computation, for $\sigma = P, Q$, we can see that $\theta_{F_n}(\sigma)=0$ since $\rho_3(\sigma) s(\rho_2(\sigma^{-1}))$ satisfies
\[\begin{split}
   [s_{ij}(x_l)]_3 \mapsto [s_{ij}(x_l)]_3
  \end{split}\]
for any $1 \leq i, j \leq m$ and $1 \leq l \leq n$.
Consider $\rho_3(S) s(\rho_2(S^{-1}))$. By using (\ref{eq-pooh}), for any $1 \leq i, j \leq m$, we see
\[\begin{split}
  s(\rho_2(S^{-1}))([s_{ij}(x_1)]_3) & = - [s_{ij}(x_1)]_3, \\
  \rho_3(S)(-[s_{ij}(x_1)]_3) & =- [s_{ij}(x_1^{-1})]_3 = [s_{ij}(x_1)]_3 - \sum_{k=1}^m [s_{ik}(x_1)s_{kj}(x_1)]_3.
  \end{split}\]
Hence we obtain
\[\begin{split}
   \theta_{F_n}(S) &= - \sum_{(i,j) \neq (m,m)} \sum_{k=1}^m s_{ij}(x_1)^* \otimes [s_{ik}(x_1)s_{kj}(x_1)]_3.
  \end{split}\]
where $s_{ij}(x_l)^*$ is the dual basis  of $s_{ij}(x_l)$ in $\mathrm{Hom}_{\Q}(\mathrm{gr}^1(\mathfrak{J}_{F_n}),\Q)$.
Similarly, we can obtain
\[\begin{split}
   \theta_{F_n}(U) &= - \sum_{(i,j) \neq (m,m)} \sum_{k=1}^m s_{ij}(x_1)^* \otimes ([s_{ik}(x_2)s_{kj}(x_2)
      + [s_{ik}(x_1)s_{kj}(x_2)]_3).
  \end{split}\]

\vspace{0.5em}

In the following, we construct the crossed homomorphism $f_1 : \mathrm{Aut}\,F_n \rightarrow H_{\Q}^* \otimes_{\Q} \Lambda^2 H_{\Q}$ and
$f_2 : \mathrm{Aut}\,F_n \rightarrow H_{\Q}$ from $\theta_{F_n}$, and study relations between them and $f_K$ and $f_M$.
The images of Nielsen's generators of $\mathrm{Aut}\,F_n$ by $f_K$ and $f_M$ are given by
\[ f_K(\sigma) = \begin{cases}
                    -x_1^* \otimes x_1 \wedge x_2, \hspace{0.5em} & \sigma=U, \\
                    0, \hspace{0.5em} & \sigma=P, Q, S,
                 \end{cases} \hspace{1em}
 f_M(\sigma) = \begin{cases}
                    -x_1, \hspace{0.5em} & \sigma=S, \\
                    0, \hspace{0.5em} & \sigma=P, Q, U.
                 \end{cases}\]
(For details, see \cite{S01} and \cite{S22}.)
Recall that
\[ \mathrm{gr}^1(\mathfrak{J}_{F_n}) \cong H_{\Q}^{\oplus (m^2-1)}, 
  \hspace{1em} \mathrm{gr}^2(\mathfrak{J}_{F_n}) \cong (S^2 H_{\Q})^{\oplus (m^2-1)} \oplus (H_{\Q}^{\otimes 2})^{\oplus \binom{m^2-1}{2}}, \]
and that $T_1=\{ s_{ij}(x_l) \,|\, (i,j) \neq (m,m), \,\,\, 1 \leq l \leq n \}$ and $T_2$ are basis of them respectively.
Let $p_1 : \mathrm{Hom}_{\Q}(\mathrm{gr}^1(\mathfrak{J}_{F_n}), \mathrm{gr}^2(\mathfrak{J}_{F_n})) \rightarrow \mathrm{Hom}_{\Q}(H_{\Q}, \mathrm{gr}^2(\mathfrak{J}_{F_n}))$
be the homomorphism induced from the inclusion map $H_{\Q} \rightarrow \mathrm{gr}^1(\mathfrak{J}_{F_n})$ defined by
\[ \sum_{l=1}^n c_{l} x_{l} \mapsto \sum_{l=1}^n c_{l} s_{11}(x_l). \]
Let $p_2 : \mathrm{Hom}_{\Q}(H_{\Q}, \mathrm{gr}^2(\mathfrak{J}_{F_n})) \rightarrow \mathrm{Hom}_{\Q}(H_{\Q}, H_{\Q}^{\otimes 2})$
be the homomorphism induced from the projection $\mathrm{gr}^2(\mathfrak{J}_{F_n}) \rightarrow H_{\Q}^{\otimes 2}$ defined by
\[ \sum_{t \in T_2} c_t t \mapsto \sum_{1 \leq i \leq j \leq n} c_{s_{12}(x_i) s_{21}(x_j)} x_i \otimes x_j. \]
Let $p_3 : \mathrm{Hom}_{\Q}(H_{\Q}, H_{\Q}^{\otimes 2}) \rightarrow \mathrm{Hom}_{\Q}(H_{\Q}, \Lambda^2 H_{\Q})$ be the homomorphism
induced from the homomorphism $H_{\Q}^{\otimes 2} \rightarrow \Lambda^2 H_{\Q}$ defined by
\[ x_i \otimes x_j \mapsto x_i \wedge x_j \]
for any $1 \leq i < j \leq n$.
Then the composition map $p_3 \circ p_2 \circ p_1$, denoted by $p$, is an $\mathrm{Aut}\,F_n$-equivariant homomorphism.
Set
\[ f_1 := p \circ \theta_{F_n} : \mathrm{Aut}\,F_n \rightarrow H_{\Q}^* \otimes_{\Q} \Lambda^2 H_{\Q}. \]

\vspace{0.5em}

On the other hand, 
Let $q_2 : \mathrm{Hom}_{\Q}(H_{\Q}, \mathrm{gr}^2(\mathfrak{J}_{F_n})) \rightarrow \mathrm{Hom}_{\Q}(H_{\Q}, S^2 H_{\Q})$
be the homomorphism induced from the projection $\mathrm{gr}^2(\mathfrak{J}_{F_n}) \rightarrow S^2 H_{\Q}$ defined by
\[ \sum_{t \in T_2} c_t t \mapsto \sum_{1 \leq i \leq j \leq n} c_{s_{11}(x_i)s_{11}(x_j)} x_i x_j. \]
Let $q_3 : \mathrm{Hom}_{\Q}(H_{\Q}, S^2 H_{\Q}) \rightarrow \mathrm{Hom}_{\Q}(H_{\Q}, H_{\Q}^{\otimes 2})$ be the homomorphism
induced from the homomorphism $S^2 H_{\Q} \rightarrow H_{\Q}^{\otimes 2}$ defined by
\[ x_i x_j \mapsto x_i \otimes x_j + x_j \otimes x_i \]
for any $1 \leq i \leq j \leq n$. Let $q_4 : \mathrm{Hom}_{\Q}(H_{\Q}, H_{\Q}^{\otimes 2}) \cong H_{\Q}^* \otimes_{\Q} H_{\Q}^{\otimes 2} \rightarrow H_{\Q}$
be the contraction map with respect to the first and the second component.
Then the composition map $q_4 \circ q_3 \circ q_2 \circ p_1$, denoted by $q$, is an $\mathrm{Aut}\,F_n$-equivariant homomorphism.
Set
\[ f_2 := q \circ \theta_{F_n} : \mathrm{Aut}\,F_n \rightarrow H_{\Q}. \]
Finally, set
\[ x := x_1 + x_2 + \cdots + x_n \in H_{\Q}, \]
and let $\delta_x$ be the principal crossed homomorphism associated to $x$.
Then we obtain the following.
\begin{thm}\label{T-Haruhi}
For any $n \geq 2$,
\[\begin{split}
 f_K = f_1, \,\,\, f_M = - f_2 + \delta_x
  \end{split}\]
as crossed homomorphisms. Here $\delta_x$ is the principle derivation associated to $x \in H_{\Q}$.
\end{thm}
\textit{Proof.}
Since
\[ f_1(\sigma) = \begin{cases}
                    -x_1^* \otimes x_1 \wedge x_2, \hspace{0.5em} & \sigma=U, \\
                    0, \hspace{0.5em} & \sigma=P, Q, S,
                 \end{cases} \hspace{1em}
 f_2(\sigma) = \begin{cases}
                    0, \hspace{1em} & \sigma=P, Q, \\
                    -x_1, \hspace{1em} & \sigma=S, \\
                    -x_2, \hspace{1em} & \sigma=U,
                 \end{cases}\]
we can obtain the required results by direct computation. $\square$

\vspace{0.5em}

Since $f_K$ and $f_M$ define non-trivial cohomology classes in $H^1(\mathrm{Aut}\,F_n, H_{\Q}^* \otimes_{\Q} \Lambda H_{\Q})= \Q^{\otimes 2}$ and
$H^1(\mathrm{Aut}\,F_n, H_{\Q})= \Q$ by \cite{S01} and \cite{S10} respectively,
we see that $\theta_{F_n}$ defines non-trivial cohomology class in
$H^1(\mathrm{Aut}\,F_n, \mathrm{Hom}_{\Q}(\mathrm{gr}^1(\mathfrak{J}_{F_n}), \mathrm{gr}^2(\mathfrak{J}_{F_n})))$.
We remark that
in our forthcoming paper \cite{S23}, we computed
\[ H^1(\mathrm{Aut}\,F_n, H_{\Q}^* \otimes_{\Q} S^2 H_{\Q}) \cong \Q. \]
Thus we have
\[\begin{split}
  H^1(& \mathrm{Aut}\,F_n, \mathrm{Hom}_{\Q}(\mathrm{gr}^1(\mathfrak{J}_{F_n}), \mathrm{gr}^2(\mathfrak{J}_{F_n}))) \\
  & \cong H^1(\mathrm{Aut}\,F_n, H_{\Q}^* \otimes_{\Q} S^2 H_{\Q})^{\oplus (m^2-1)^2} \\
  & \hspace{4em} \bigoplus H^1(\mathrm{Aut}\,F_n, H_{\Q}^* \otimes_{\Q} H_{\Q}^{\otimes 2})^{\oplus (m^2-1) \binom{m^2-1}{2} } \\
  & \cong \Q^{\oplus \frac{1}{2}(m^2-1)^2(3m^2-4)}.
  \end{split}\]

\section{On the $\mathrm{SL}(m,\C)$-representation algebra of $H$}\label{S-FAb}

\vspace{0.5em}

In this section, we study the $\mathrm{SL}(m,\C)$-representation algebra for free abelian groups.
Recall that $H$ is the free abelian group of rank $n$ with basis $x_1, \ldots, x_n$.
It seems to be quite difficult to obtain a basis of $\mathrm{gr}^k(\mathfrak{J}_{H})$ for a general $k \geq 1$.
Here, we give basis of $\mathrm{gr}^k(\mathfrak{J}_{H})$ for $1 \leq k \leq 2$. Then, by using it, we study the extension of $\eta_1$.

\vspace{0.5em}

By Lemma {\rmfamily \ref{L-Dream}}, we have the surjective homomorphism
\[ \overline{\mathfrak{a}} : \mathfrak{R}_{\Q}^m(F_n) \rightarrow \mathfrak{R}_{\Q}^m(H) \]
induced from the abelianization $\mathfrak{a} : F_n \rightarrow H$. The above map $\overline{\mathfrak{a}}$ also induces surjective homomorphisms
\[ \overline{\mathfrak{a}}^k : \mathrm{gr}^k(\mathfrak{J}_{F_n}) \rightarrow \mathrm{gr}^k(\mathfrak{J}_{H}) \]
for any $k \geq 1$. We study the kernel of $\overline{\mathfrak{a}}^k$ for $1 \leq k \leq 2$.

\vspace{0.5em}

Denote by $F_{ij}(z)$ the $m \times m$-matrix whose $(i,j)$-entry is $z \in \C$ and the other entries zero.
Furthermore, we denote by $E_i(z) \in \mathrm{SL}(m,\C)$ the matrix whose $(i,i)$-entry is $z \in \C \setminus \{ 0 \}$, the other diagonal entries one, and
the other entries zero.

\subsection{The graded quotients $\mathrm{gr}^1(\mathfrak{J}_{H})$}\label{Ss-FAb-1}
\hspace*{\fill}\ 

\vspace{0.5em}

First, we show that $\overline{\mathfrak{a}}^1$ is isomorphism.
From Propositions {\rmfamily \ref{P-dim-2}} and {\rmfamily \ref{P-Rythm}},
we have $\mathrm{gr}^1(\mathfrak{J}_{F_n}) \cong H_{\Q}^{\oplus m^2-1}$ and
\[ T_1 = \{ s_{ij}(x_l) \,|\, 1 \leq i, j \leq m, \,\,\, (i,j) \neq (m,m), \,\,\, 1 \leq l \leq n \} \pmod{\mathfrak{J}_{F_n}^2} \]
forms a basis of $\mathrm{gr}^1(\mathfrak{J}_{F_n})$. 
Let $\overline{x_1}, \ldots, \overline{x_n}$ be the coset classes of the basis $x_1, \ldots, x_n$ of $F_n$.
Then, for each $1 \leq i, j \leq m$ and $1 \leq l \leq n$, we have
$s_{ij}(\overline{x_l}) = s_{ij}(\mathfrak{a}(x_l)) = \overline{\mathfrak{a}}(s_{ij}(x_l))$. Set
\[ \overline{T_1} := \{ s_{ij}(\overline{x_l}) \,|\, 1 \leq i, j \leq m, \,\,\, (i,j) \neq (m,m), \,\,\, 1 \leq l \leq n \} \subset \mathfrak{R}_{\Q}^m(H). \]
\begin{pro}\label{P-Jintan}
The set $\overline{T_1} \pmod{\mathfrak{J}_{H}^2}$
forms a basis of $\mathrm{gr}^1(\mathfrak{J}_{H})$ as a $\Q$-vector space.
\end{pro}
\textit{Proof.}
It suffices to show that $\overline{T_1} \pmod{\mathfrak{J}_{H}^2}$ is linearly independent over $\Q$. Set
\[ f := \sum_{(i,j) \neq (m,m)} \sum_{1 \leq l \leq n} c_{ij, l} s_{ij}(\overline{x_l}) \in \mathfrak{R}_{\Q}^m(H) \]
for $c_{ij,l} \in \Q$, and assume that $f \equiv 0 \pmod{\mathfrak{J}_{H}^2}$. Take any $(i,j) \neq (m,m)$, and fix it.

\vspace{0.5em}

If $i \neq j$, take any $z_l \in \C$ for $1 \leq l \leq n$, and
consider the representation $\rho_1 : H \rightarrow \mathrm{SL}(m,\C)$ defined by
\[ x_l \mapsto E_n + F_{ij}(z_l). \]
Then we see that
\[ f(\rho) = \sum_{1 \leq l \leq n} c_{ij,l} z_l. \]
The assumption $f \equiv 0 \pmod{\mathfrak{J}_{H}^2}$ means
\[ \sum_{1 \leq l \leq n} c_{ij,l} z_l = \text{terms of degree $\geq 2$}. \]
Since we can take $z_l \in \C$ arbitrary, we see that $c_{ij,l}=0$ for any $1 \leq l \leq n$.

\vspace{0.5em}

If $i=j$, take any $z_l \in D := \{ z \in \C \,|\, |z|<1 \}$ for $1 \leq l \leq n$, and
consider the representation $\rho_2 : H \rightarrow \mathrm{SL}(m,\C)$ defined by
\[ x_l \mapsto E_{i}(1+z_l) E_m((1+z_l)^{-1}). \]
Then by the same argument as above, we see $c_{ii,l}=0$ for any $1 \leq l \leq n$.
Thus we obtain the required result. $\square$

\subsection{The graded quotients $\mathrm{gr}^2(\mathfrak{J}_{H})$}\label{Ss-FAb-2}
\hspace*{\fill}\ 

\vspace{0.5em}

Next we consider the case of $\mathrm{gr}^2(\mathfrak{J}_{H})$. The basic idea of the strategy to find a basis of $\mathrm{gr}^2(\mathfrak{J}_{H})$
is the same as that of $\mathrm{gr}^1(\mathfrak{J}_{H})$.
From Propositions {\rmfamily \ref{P-dim-2}} and {\rmfamily \ref{P-Rythm}}, we have
\[ \mathrm{gr}^2(\mathfrak{J}_{F_n}) \cong (S^2 H_{\Q})^{\oplus m^2-1} \oplus (H_{\Q}^{\otimes 2})^{\oplus \binom{m^2-1}{2}} \]
and,
\[\begin{split}
  T_2 & = \{ s_{ij}(x_{p})s_{ij}(x_{q}) \,|\, 1 \leq i, j \leq m, \,\,\, (i,j) \neq (m,m), \,\,\, 1 \leq p \leq q \leq n \} \\
      & \hspace{1em} \cup \{ s_{ij}(x_{p})s_{hk}(x_{q}) \,|\, (i,j) <_{\mathrm{lex}} (h,k), \,\,\, (i,j), (h,k) \neq (m,m), \,\,\, 1 \leq p, q \leq n \}
  \end{split}\]
$\pmod{\mathfrak{J}_{F_n}^3}$ forms a basis of $\mathrm{gr}^2(\mathfrak{J}_{F_n})$. We rewrite this basis according to the decomposition
$H_{\Q}^{\otimes 2} \cong \Lambda^2 H_{\Q} \oplus S^2 H_{\Q}$.
Set
\[\begin{split}
  t_{ij,hk}(p,q) & := s_{ij}(x_p) s_{hk}(x_q) - s_{ij}(x_q) s_{hk}(x_p), \\
  u_{ij,hk}(p,q) & := s_{ij}(x_p) s_{hk}(x_q) + s_{ij}(x_q) s_{hk}(x_p), \\
  v_{ij}(p,q) & : = s_{ij}(x_{p})s_{ij}(x_{q}).
 \end{split}\]
Then,
\[\begin{split}
   T_2' &:= \{ t_{ij,hk}(p,q) \,|\, (i, j) <_{\mathrm{lex}} (h, k), \,\,\, (i,j), (h,k) \neq (m,m), \,\,\, p<q \} \\
        & \hspace{2em} \cup \{ u_{ij,hk}(p,q) \,|\, (i,j) <_{\mathrm{lex}} (h, k), \,\,\, (i,j), (h,k) \neq (m,m), \,\,\, p \leq q \} \\
        & \hspace{2em} \cup \{ v_{ij}(p,q) \,|\, (i,j) \neq (m,m), \,\,\, 1 \leq p \leq q \leq n \} \pmod{\mathfrak{J}_{F_n}^3}
  \end{split}\]
also forms a basis of $\mathrm{gr}^2(\mathfrak{J}_{F_n})$.

\vspace{0.5em}

Now, we observe some elements in the kernel of $\overline{\mathfrak{a}}^2$.
For any $1 \leq i, j \leq m$, and any $1 \leq p, q \leq n$, the equation
$\overline{\mathfrak{a}}(s_{ij}(x_p x_q)) = \overline{\mathfrak{a}}(s_{ij}(x_q x_p))$ induces the fact that
\[ R_{ij}(p,q) := \sum_{k=1}^m t_{ik, kj}(p,q) \pmod{\mathfrak{J}_{F_n}^3} \in \mathrm{Ker}(\overline{\mathfrak{a}}^2). \]
More precisely, if we describe the above by using the elements of $T_2'$,
\[\begin{split}
   R_{ij}(p,q) & = -t_{1j, i1}(p,q) - \cdots - t_{i-1 j, i i-1}(p,q) + t_{ii, ij}(p,q) + \cdots + t_{im, mj}(p,q), \hspace{1em} (i \neq 1), \\
   R_{1j}(p,q) & = t_{11, 1j}(p,q) + \cdots + t_{1m, mj}(p,q), \hspace{1em} (j \neq 1).
  \end{split}\]
Set
\[\begin{split}
   \overline{t}_{ij,hk}(p,q) & := \overline{\mathfrak{a}}(t_{ij, hk}(p,q)), \hspace{1em} \overline{u}_{ij,hk}(p,q) := \overline{\mathfrak{a}}(u_{ij, hk}(p,q)), \\
   \overline{v}_{ij}(p,q) & := \overline{\mathfrak{a}}(v_{ij}(p,q)).
  \end{split}\]
From the above observation, we see that we can remove $\overline{t}_{1j,i1}(p,q)$, $\overline{t}_{11,i1}(p,q)$ and $\overline{t}_{11,1j}(p,q)$ for $i, j \neq 1$
from the generating set $\mathfrak{a}^2(T_2')$ of $\mathrm{gr}^2(\mathfrak{J}_{H})$ by using $R_{ij}(p,q)$ and $R_{1j}(p,q)$ respectively.
Define the index sets $I, J$ by
\[\begin{split}
  I & := \{ (i,j,h,k) \,|\, (i,j) <_{\mathrm{lex}} (h,k), \,\,\, (i,j), (h,k) \neq (m,m) \}, \\
  J & := I \setminus \{ (1, j, i, 1), (1,1,i,1), (1,1,1,j) \,|\, i, j \neq 1 \}.
  \end{split}\]
Set
\[\begin{split}
  Y &:= \{ \overline{t}_{ij,hk}(p,q) \,|\, (i,j,h,k) \in J, \,\,\, p<q \} \cup \{ \overline{u}_{ij,hk}(p,q) \,|\, (i,j,h,k) \in I, \,\,\, p \leq q \} \\
        & \hspace{2em} \cup \{ \overline{v}_{ij}(p,q) \,|\, (i,j) \neq (m,m), \,\,\, p \leq q \}.
  \end{split}\]

\begin{pro}\label{P-haruhi}
The set $Y \pmod{\mathfrak{J}_{H}^3}$ forms a basis of $\mathrm{gr}^2(\mathfrak{J}_{H})$ as a $\Q$-vector space.
\end{pro}
\textit{Proof.}
It suffices to show that $Y \pmod{\mathfrak{J}_{H}^3}$ is linearly independent over $\Q$.
Set
\[\begin{split}
   f & := \sum_{(i,j,h,k) \in J} \sum_{p<q} \alpha_{ij,hk}(p,q) \overline{t}_{ij,hk}(p,q) + \sum_{(i,j,h,k) \in I} \sum_{p \leq q} \beta_{ij,hk}(p,q) \overline{u}_{ij,hk}(p,q) \\
     & \hspace{2em} + \sum_{(i,j) \neq (m,m)} \sum_{p \leq q} \gamma_{ij}(p,q) \overline{v}_{ij}(p,q) \in \mathfrak{R}_{\Q}^m(H)
  \end{split}\]
for $\alpha_{ij,hk}(p,q), \beta_{ij,hk}(p,q), \gamma_{ij}(p,q) \in \Q$, and assume $f \equiv 0 $.

\vspace{0.5em}

{\bf Step 1.} The proof of $\gamma_{ij}(p,q)=0$.

\vspace{0.5em}

Take any $(i,j) \neq (m,m)$, and fix it. If $i \neq j$,
by using the representation $\rho_1 : H \rightarrow \mathrm{SL}(m,\C)$ defined in the proof of Proposition {\rmfamily \ref{P-Jintan}},
and by using $f \equiv 0 \pmod{\pmod{\mathfrak{J}_{H}^3}}$,
we see that
\[ \sum_{p \leq q} \gamma_{ij}(p,q) z_p z_q = \text{terms of degree $\geq 3$}. \]
Since we can take $z_l \in \C$ arbitrary, we see that $\gamma_{ij}(p,q)=0$ for any $1 \leq p \leq q \leq n$.
Similarly, if $i=j$, by considering the representation $\rho_2 : H \rightarrow \mathrm{SL}(m,\C)$
defined in the proof of Proposition {\rmfamily \ref{P-Jintan}}, we can obtain
$\gamma_{ii}(p,q)=0$ for any $1 \leq p \leq q \leq n$.

\vspace{0.5em}

{\bf Step 2.} The proof of $\beta_{ij}(p,q)=0$.

\vspace{0.5em}

For any $1 \leq l \leq n$, take any $z_{l1}, \ldots, z_{l m-1} \in \C \setminus \{ 0 \}$, and set $z_{lm} := (z_{l1} \cdots z_{l m-1})^{-1}$.
For any $B=(b_{ij}) \in \mathrm{GL}(m,\C)$, consider the representation
$\rho_3 : H \rightarrow \mathrm{SL}(m,\C)$ defined by
\[ x_l \mapsto \dfrac{1}{|B|} \Bigg{[} B \begin{pmatrix} z_{l1} &        &      \\
                                     & \ddots &      \\
                                     &        & z_{l m} \end{pmatrix} \widetilde{B} \Bigg{]} \]
for any $1 \leq l \leq n$ where $\widetilde{B}=(\tilde{b}_{ji})$ is the adjugate matrix of $B$. By considering the power series expansion
among $z_{lr}-1$, we can see that for any $1 \leq l \leq n$, the $(i,j)$-entry of $\rho_3(x_l)-E_m$ is given by
\[ \dfrac{1}{|B|} \Big{(} \sum_{r=1}^{m-1} (b_{ir} \tilde{b}_{jr} - b_{im} \tilde{b}_{jm}) (z_{lr} -1) + (\text{terms of degree $\geq 2$}) \Big{)}. \]
Hence we have
\[\begin{split}
  f(\rho_3) & = \dfrac{1}{|B|^2} \Bigg{\{} \sum_{(i,j,h,k) \in J} \sum_{p<q} \alpha_{ij,hk}(p,q) \sum_{r,s=1}^{m-1}
                \Big{(} (b_{ir} \tilde{b}_{jr} - b_{im} \tilde{b}_{jm})(b_{hs} \tilde{b}_{ks} - b_{hm} \tilde{b}_{km}) \\
            & \hspace{10em} - (b_{is} \tilde{b}_{js} - b_{im} \tilde{b}_{jm})(b_{hr} \tilde{b}_{kr} - b_{hm} \tilde{b}_{km}) \Big{)} (z_{pr}-1)(z_{qs}-1) \\
            & \hspace{1em} + \sum_{(i,j,h,k) \in I} \sum_{p \leq q} \beta_{ij,hk}(p,q) \sum_{r,s=1}^{m-1}
                \Big{(} (b_{ir} \tilde{b}_{jr} - b_{im} \tilde{b}_{jm})(b_{hs} \tilde{b}_{ks} - b_{hm} \tilde{b}_{km}) \\
            & \hspace{10em} + (b_{is} \tilde{b}_{js} - b_{im} \tilde{b}_{jm})(b_{hr} \tilde{b}_{kr} - b_{hm} \tilde{b}_{km}) \Big{)} (z_{pr}-1)(z_{qs}-1) 
              \Bigg{\}} \\
            & \hspace{1em} + (\text{terms of degree $\geq 3$}).
  \end{split}\]
Since $f \equiv 0 \pmod{\pmod{\mathfrak{J}_{H}^3}}$, we see that each coefficient of $(z_{pr}-1)(z_{qs}-1)$ is equal to zero.
In particular, for any $p \leq q$, by observing the coefficient of $(z_{p1}-1)(z_{q1}-1)$, we have
\[ \sum_{(i,j,h,k) \in I} 2 \beta_{ij,hk}(p,q) (b_{i1} \tilde{b}_{j1} - b_{im} \tilde{b}_{jm})(b_{h1} \tilde{b}_{k1} - b_{hm} \tilde{b}_{km}) = 0. \]
Now we consider the linear ordering on the monomials among $b_{ij}$ induced from the lexicographic ordering on $\{ (i,j) \,|\, 1 \leq i,j \leq m \}$.
Then the minimum term in $(b_{i1} \tilde{b}_{j1} - b_{im} \tilde{b}_{jm})$ is
\[ b_{im}b_{11} b_{22} \cdots b_{j-1 j-1} b_{j+1 j} b_{j+2 j+1} \cdots b_{m m-1} \]
coming from $b_{im}\tilde{b}_{jm}$, and that in $(b_{h1} \tilde{b}_{k1} - b_{hm} \tilde{b}_{km})$ is
\[ b_{hm}b_{11} b_{22} \cdots b_{k-1 k-1} b_{k+1 k} b_{k+2 k+1} \cdots b_{m m-1} \]
coming from $b_{hm}\tilde{b}_{km}$. We focus on the term
\[ b_{im}b_{11} \cdots b_{j-1 j-1} b_{j+1 j} \cdots b_{m m-1} \cdot b_{hm}b_{11} \cdots b_{k-1 k-1} b_{k+1 k} \cdots b_{m m-1}. \]

\vspace{0.5em}

Consider $(i,j,h,k)=(1,1,1,2)$. In $f(\rho_3)$, the coefficient of
\[ b_{1m} b_{21}b_{32} \cdots b_{m m-1} \cdot b_{1m}b_{11} \cdots b_{32} \cdots b_{m m-1} \]
is $\beta_{11,12}(p,q)$, and hence $\beta_{11,12}(p,q)=0$ since we can take $B=(b_{ij}) \in \mathrm{GL}(n,\C)$ arbitrary.
Next, consider $(i,j,h,k)=(1,1,1,3)$.
By observing the coefficient of
\[ b_{1m}b_{11} b_{21} \cdots b_{m m-1} \cdot b_{1m}b_{11} b_{22} b_{43} \cdots b_{m m-1}, \]
we see $\beta_{11,13}(p,q)=0$. Similarly, by an inductive argument, we see
\[ \beta_{11,hk}(p,q)=0, \,\,\, \beta_{12,hk}(p,q)=0, \,\, \ldots, \,\, \beta_{m m-2, m m-1}=0. \]

\vspace{0.5em}

{\bf Step 3.} The proof of $\alpha_{ij}(p,q)=0$.

\vspace{0.5em}

Take any $(i,j,h,k) \in J$, and set $N:= \sharp \{ i, j, h, k \}$.

\vspace{0.5em}

{\bf (1)} The case of $N=4$.

\vspace{0.5em}

For any $1 \leq l \leq n$, take any $z_l, w_l \in \C$, and
consider the representation $\rho_4 : H \rightarrow \mathrm{SL}(m,\C)$ defined by
\[ x_l \mapsto E_m + F_{ij}(z_l) + F_{hk}(w_l). \]
Then $f(\rho_4)=\alpha_{ij,hk}(p,q)(z_pw_q-z_qw_p)$.
By an argument similar to the above, we see $\alpha_{ij,hk}(p,q)=0$ for any $p \leq q$.

\vspace{0.5em}

{\bf (2)} The case of $N \leq 3$.

\vspace{0.5em}

We have to consider the following cases of indices:

\vspace{0.5em}

\begin{center}
{\renewcommand{\arraystretch}{1.3}
\begin{tabular}{|c|c||c|c||c|c|} \hline
  {\bf (i)}   & $(i,i,j,h)$ & {\bf (ii)}    & $(i,j,i,h)$ & {\bf (iii)}  & $(i,j,h,i)$ \\ \hline
  {\bf (iv)}  & $(i,j,j,h)$ & {\bf (v)}     & $(i,j,h,j)$ & {\bf (vi)}   & $(i,j,h,h)$ \\ \hline \hline
  {\bf (vii)} & $(i,i,i,j)$ & {\bf (viii)}  & $(i,i,j,i)$ & {\bf (ix)}   & $(i,j,i,i)$ \\ \hline
  {\bf (x)}   & $(i,j,j,j)$ & {\bf (xi)}    & $(i,i,j,j)$ & {\bf (xii)}  & $(i,j,j,i)$ \\ \hline
\end{tabular}}
\end{center}

\vspace{0.5em}

{\bf (ii) and (v).} For any $1 \leq l \leq n$, take any $z_l, w_l \in \C$.
In the case where $(i,j,i,h) \in J$, 
define the representation $\rho_5 : H \rightarrow \mathrm{SL}(m,\C)$ by
\[ x_l \mapsto E_m + F_{ij}(z_l) + F_{ih}(w_l). \]
Then we can see $\alpha_{ij,ih}(p,q)=0$ for any $p \leq q$.
Similarly, in the case where $(i,j,h,j) \in J$, from the representation $\rho_6 : H \rightarrow \mathrm{SL}(m,\C)$ defined by
\[ x_l \mapsto E_m + F_{ij}(z_l) + F_{hj}(w_l), \]
we can obtain $\alpha_{ij,hj}(p,q)=0$ for any $p \leq q$.

\vspace{0.5em}

Next, we consider the other cases.
For any $1 \leq l \leq n$, take any $z_l, w_l \in D$.
For any $b_{12}, b_{13}, b_{21}, b_{23}, b_{31}, b_{32} \in \C$ such that $\Delta := b_{12}b_{23}b_{31}+ b_{13}b_{32}b_{21} \neq 0$, define the $3 \times 3$ matrix $C_l$ by
{\footnotesize
\[\begin{split}
   C_l & := \begin{pmatrix} 0 & b_{12} & b_{13} \\ b_{21} & 0 & b_{23} \\ b_{31} & b_{32} & 0 \end{pmatrix}
        \begin{pmatrix} z_l & 0 & 0 \\ 0 & w_l & 0 \\ 0 & 0 & (z_l w_l)^{-1} \end{pmatrix}
        \begin{pmatrix} 0 & b_{12} & b_{13} \\ b_{21} & 0 & b_{23} \\ b_{31} & b_{32} & 0 \end{pmatrix}^{-1} \\
     & = \dfrac{1}{\Delta} \begin{pmatrix} b_{12}b_{23}b_{31} w_l+ b_{13}b_{32}b_{21}(z_l w_l)^{-1} & -b_{12}b_{13}b_{31} (w_l-(z_lw_l)^{-1})
                                           & b_{12}b_{13}b_{21}(w_l-(z_lw_l)^{-1}) \\
         b_{21}b_{23}b_{32}(-z_l+(z_lw_l)^{-1}) & b_{12}b_{23}b_{31}(z_lw_l)^{-1} + b_{13}b_{32}b_{21} z_l & b_{12}b_{21}b_{23}(z_l-(z_lw_l)^{-1}) \\
         b_{23} b_{31} b_{32}(-z_l+w_l) & b_{13}b_{31}b_{32}(z_l-w_l) & b_{12}b_{23}b_{31} z_l + b_{13}b_{32}b_{21} w_l \end{pmatrix}.
  \end{split}\]
}
By observing the power series expansion of each entry of $C_l-E_3$ with respect to $z_l' := z_l-1$ and $w_l' := w_l-1$, and by writing down the degree one part
of each entry, we have
{\footnotesize
\[ \dfrac{1}{\Delta} \begin{pmatrix} b_{12}b_{23}b_{31} w_l' - b_{13}b_{32}b_{21}(z_l'+ w_l') & -b_{12}b_{13}b_{31} (z_l'+2w_l')
                                           & b_{12}b_{13}b_{21}(z_l' +2 w_l') \\
         -b_{21}b_{23}b_{32}(w_l'+2z_l') & -b_{12}b_{23}b_{31}(z_l' + w_l') + b_{13}b_{32}b_{21} z_l' & b_{12}b_{21}b_{23}(2z_l'+w_l') \\
         b_{23} b_{31} b_{32}(-z_l'+w_l') & b_{13}b_{31}b_{32}(z_l'-w_l') & b_{12}b_{23}b_{31} z_l' + b_{13}b_{32}b_{21} w_l' \end{pmatrix}. \]
}

\vspace{0.5em}

Take any $1 \leq i < j< k \leq m$, and fix it. For any $C=(c_{ij}) \in \mathrm{SL}(3,\C)$, set
\[ \begin{split} \widehat{C} := 
             \bordermatrix{
                           &       &\underline{i}    &\underline{j} &\underline{k} &        \\
                           &E      & O               &\cdots        & \cdots       &O       \\
             \underline{i} &O      &c_{11}           & c_{12}       & c_{13}       &\vdots  \\
             \underline{j} &\vdots &c_{21}           & c_{22}       & c_{23}       &\vdots  \\
             \underline{k} &\vdots &c_{31}           & c_{32}       & c_{33}       & O      \\
                           &O      &\cdots           &\cdots        &  O           & E       } \in \mathrm{SL}(m,\C)
 \end{split} \]
where $E$ denotes the identity matrix.
Namely, $\widehat{C}$ is obtained from the identity matrix with $(i,i)$, $(i,j)$, $(i,k), \ldots$ entries replaced with $c_{11}$, $c_{12}$, $c_{13}, \ldots$ respectively.
Define the representation $\rho_7 = \rho_7(i,j,k) : H \rightarrow \mathrm{SL}(m,\C)$ by $x_l \mapsto \widehat{C_l}$ for any $1 \leq l \leq n$.
Then by considering the power series expansion with respect to $z_l'$ and $w_l'$, we see that
$f(\rho_7)$ can be written as a power series of $z_l'$ and $w_l'$ such that the minimal degree of the monomials of $f(\rho_7)$ is greater than or equal to $2$.
From the assumption $f \equiv 0 \pmod{\pmod{\mathfrak{J}_{H}^3}}$, we can see that each of coefficients of the monomials of degree $2$ is equal to zero.

\vspace{0.5em}

{\bf (xi).} Consider the case $(i,j,h,k)=(i,i,j,j)$.
Recall that $(i,i), (j,j) \neq (m,m)$.
For any $p \leq q$, by observing the coefficient of $b_{12}^2b_{23}^2b_{31}^2 z_p' w_q'$ of $f(\rho_7(i,j,m))$, we see
\[ \alpha_{ii,jj}(p,q) - \alpha_{ii,mm}(p,q) + \alpha_{jj,mm}(p,q) = 0. \]
Thus we obtain $\alpha_{ii,jj}(p,q)=0$.

\vspace{0.5em}

{\bf (xii).} Consider the case $(i,j,h,k)=(i,j,j,i)$. Recall that $i \neq 1$.
For any $p \leq q$, by observing the coefficient of $b_{12} b_{13}b_{21} b_{23} b_{31}^2 z_p' w_q'$ of $f(\rho_7(1,i,j))$, we see
\[ - \alpha_{1i,11}(p,q) + \alpha_{1j,j1}(p,q) - \alpha_{ij,ji}(p,q) = 0, \]
and hence $\alpha_{ij,ji}(p,q) = 0$ since $\alpha_{1i,11}(p,q) = \alpha_{1j,j1}(p,q) = 0$.

\vspace{0.5em}

{\bf (vi) and (i).} In general,
for any $p \leq q$, from the coefficients of $b_{12}^2 b_{13} b_{23} b_{31}^2 z_p' w_q'$, $b_{12} b_{13}^2 b_{21}^2 b_{32} z_p' w_q'$ and
$b_{13} b_{21}^2 b_{23} b_{32}^2 z_p' w_q'$ of $f(\rho_7(i,j,k))$, we obtain
\begin{eqnarray}
 \alpha_{ii,ij}(p,q) - \alpha_{ij,jj}(p,q) +2 \alpha_{ij,kk}(p,q) = 0 \label{eq-haruhi-1} \\
 - \alpha_{ii,ik}(p,q) + \alpha_{ik,kk}(p,q) +2 \alpha_{ik,jj}(p,q) = 0 \label{eq-haruhi-2} \\
 -\alpha_{ii,ji}(p,q) + \alpha_{ji,jj}(p,q) -2 \alpha_{ji,kk}(p,q) = 0 \label{eq-haruhi-3}
\end{eqnarray}
respectively.
By considering the case where $k=m$ in (\ref{eq-haruhi-1}), we have
\begin{equation}
 \alpha_{ii,ij}(p,q) = \alpha_{ij,jj}(p,q), \label{eq-haruhi-4}
\end{equation}
and from (\ref{eq-haruhi-1}) again $\alpha_{ij,kk}(p,q) = 0$. 
By the same argument as above, from (\ref{eq-haruhi-3}) we obtain
\begin{equation}
 \alpha_{ji,kk}(p,q)=0, \hspace{1em} \alpha_{ii,ji}(p,q) = \alpha_{ji,jj}(p,q), \label{eq-haruhi-5}
\end{equation}
On the other hand, from (\ref{eq-haruhi-4}) and (\ref{eq-haruhi-2}), we have $\alpha_{ik,jj}(p,q)=0$.
By changing the role of indices $i, j$ and $k$ in the above argument, we can also obtain
\[ \alpha_{ii,jk}(p,q) = \alpha_{ii,kj}(p,q) = \alpha_{jj,ki}(p,q) =0. \]

\vspace{0.5em}

{\bf (x) and (vii).} In general, for any $p \leq q$, from the coefficient of $b_{12} b_{13} b_{21}^2 b_{23} b_{32} z_p' w_q'$ of $f(\rho_7(i,j,k))$,
we see
\begin{equation}
 -3\alpha_{ik,ji}(p,q) + \alpha_{jj,jk}(p,q) +2 \alpha_{jk,kk}(p,q) = 0. \label{eq-haruhi-6}
\end{equation}
From this and (\ref{eq-haruhi-4}), we have
\begin{equation}
 \alpha_{ik,ji}(p,q) = \alpha_{jj,jk}(p,q)= \alpha_{jk,kk}(p,q). \label{eq-haruhi-7}
\end{equation}
Hence if $1<j<k$,
\[ \alpha_{jk,kk}(p,q) = \alpha_{1k,j1}(p,q) = 0, \]
and if $1<k$,
\[ \alpha_{1k,kk}(p,q) = \alpha_{11,1k}(p,q) = 0. \]
On the other hand, if $1<j<k$,
\[ \alpha_{jj,jk}(p,q) = \alpha_{1k,j1}(p,q) = 0. \]

\vspace{0.5em}

{\bf (viii) and (ix).} In general, for any $p \leq q$, from the coefficient of $b_{12} b_{13} b_{23} b_{31}^2 b_{32} z_p' w_q'$
of $f(\rho_7(i,j,k))$, we see
\begin{equation}
 -3\alpha_{ij,ki}(p,q) + 2 \alpha_{jj,kj}(p,q) + \alpha_{kj,kk}(p,q) = 0. \label{eq-haruhi-8}
\end{equation}
By (\ref{eq-haruhi-5}), we see
\begin{equation}
 \alpha_{ij,ki}(p,q) = \alpha_{jj,kj}(p,q)= \alpha_{kj,kk}(p,q). \label{eq-haruhi-9}
\end{equation}
Hence if $1<j<k$,
\[ \alpha_{jj,kj}(p,q)= \alpha_{kj,kk}(p,q) = \alpha_{1j,k1}(p,q)=0, \]
and if $1<k$,
\[ \alpha_{k1,kk}(p,q) = \alpha_{11,k1}(p,q) = 0. \]

\vspace{0.5em}

{\bf (iii).} For any $p \leq q$, from (\ref{eq-haruhi-7}), (\ref{eq-haruhi-9}) and the above argument, we have
$\alpha_{ij,ki}(p,q)=\alpha_{ik,ji}(p,q)=0$. Furthermore, by observing the coefficient of $b_{13} b_{21} b_{23} b_{31} b_{32}^2 z_p' w_q'$
of $f(\rho_7(i,j,k))$, we obtain $\alpha_{ji,kj}(p,q)=0$.

\vspace{0.5em}

{\bf (iv).} For any $p \leq q$, by a similar argument as above, from the coefficients of $b_{12}^2 b_{13} b_{21} b_{23} b_{31} z_p' w_q'$,
$b_{12} b_{13}^2 b_{21} b_{31} b_{32} z_p' w_q'$ and
$b_{12} b_{21} b_{23}^2 b_{31} b_{32} z_p' w_q'$ of $f(\rho_7(i,j,k))$, we obtain
\[ \alpha_{ij,jk}(p,q)= \alpha_{ik,kj}(p,q) = \alpha_{jk,ki}(p,q)=0 \]
respectively.

\vspace{0.5em}

This completes the proof of Proposition {\rmfamily \ref{P-haruhi}}. $\square$

\vspace{0.5em}

Remark that for each $i,j,h,k$,
$\overline{t_{ij,hk}}(p,q)$ for all $p<q$ generate $\Lambda^2 H_{\Q}$,
$\overline{u}_{ij,hk}(p,q)$ for all $p \leq q$ generate $S^2 H_{\Q}$, and
$\overline{v}_{ij}(p,q)$ for all $p \leq q$ generate $S^2 H_{\Q}$ as a $\mathrm{GL}(H_{\Q})$-module.
Thus, since $|I \setminus J|=m^2-1$,
as a corollary to Proposition {\rmfamily \ref{P-haruhi}}, we obtain the following.
\begin{pro}\label{P-haruhi-2}
For any $n \geq 2$ and $m \geq 2$, we have
\[ \mathrm{gr}^2(\mathfrak{J}_H) \cong (S^2 H_{\Q})^{\oplus \frac{1}{2}m^2(m^2-1)} \oplus (\Lambda^2 H_{\Q})^{\oplus \frac{1}{2} (m^2-1)^2 (m^2-4) }. \]
\end{pro}
We remark that this result is a generalization of our previous result in \cite{S22} for the case where $m=2$.

\subsection{The crossed homomorphism $\theta_H$}\label{Ss-H-coh}
\hspace*{\fill}\ 

\vspace{0.5em}

From the argument in Subsection {\rmfamily \ref{Ss-Ext}}, we have the crossed homomorphism
\[ \theta_H : \mathrm{Aut}\,H \rightarrow \mathrm{Hom}_{\Q}(\mathrm{gr}^1(\mathfrak{J}_H), \mathrm{gr}^2(\mathfrak{J}_H)). \]
By composing the natural map $\mathrm{Aut}\,F_n \rightarrow \mathrm{Aut}\,H$ and the above map, we obtain the crossed homomorphism
\[ \theta_H' : \mathrm{Aut}\,F_n \rightarrow \mathrm{Hom}_{\Q}(\mathrm{gr}^1(\mathfrak{J}_H), \mathrm{gr}^2(\mathfrak{J}_H)). \]
In this subsection, we study a relation between $\theta_H'$ and $f_M$.
Recall that
\[ \mathrm{gr}^1(\mathfrak{J}_H) \cong H_{\Q}^{\oplus m^2-1},
 \hspace{1em} \mathrm{gr}^2(\mathfrak{J}_H) \cong (S^2 H_{\Q})^{\oplus \frac{1}{2}m^2(m^2-1)} \oplus (\Lambda^2 H_{\Q})^{\frac{1}{2} (m^2-1)^2 (m^2-4) }, \]
and that $\{ s_{ij}(\overline{x_l}) \,|\, (i,j) \neq (m,m), \,\,\, 1 \leq l \leq n \}$ and $Y$ are basis of them as a $\Q$-vector space
respectively.

\vspace{0.5em}

By the same argument as $\theta_{F_n}$,
we can calculate the images of Nielsen's generators by $\theta_H'$.
In particular, we have
\[\begin{split}
   \theta_H(S)(s_{11}(\overline{x_1})) &= -\sum_{k=1}^m [s_{1k}(\overline{x_1})s_{k1}(\overline{x_1})]_3 \\
    & = - \overline{v}_{11}(1,1) - \sum_{k=2}^m \overline{u}_{1k,k1}(1,1) \\
   \theta_H(U)(s_{11}(\overline{x_1})) &= -\sum_{k=1}^m [s_{1k}(\overline{x_2})s_{k1}(\overline{x_2})]_3 -\sum_{k=1}^m [s_{1k}(\overline{x_1})s_{k1}(\overline{x_2})]_3 \\
    & = - \overline{v}_{11}(2,2) - \sum_{k=2}^m \overline{u}_{1k,k1}(2,2) \\
    & \hspace{2em} - \overline{v}_{11}(1,2) - \dfrac{1}{2} \sum_{k=2}^m (\overline{u}_{1k,k1}(1,2) + \overline{t}_{1k,k1}(1,2)) \\
    & = - (\overline{v}_{11}(2,2) + \overline{v}_{11}(1,2)) - \sum_{k=2}^m \Big{(} \overline{u}_{1k,k1}(2,2) + \dfrac{1}{2} \overline{u}_{1k,k1}(1,2) \Big{)}.
  \end{split}\]
Denote $s_{ij}(\overline{x_l})^*$ by the dual basis of $s_{ij}(\overline{x_l})$ in $\mathrm{Hom}_{\Q}(\mathrm{gr}^1(\mathfrak{J}_H),\Q)$.
Let $p_1' : \mathrm{Hom}_{\Q}(\mathrm{gr}^1(\mathfrak{J}_H),$ $\mathrm{gr}^2(\mathfrak{J}_H)) \rightarrow \mathrm{Hom}_{\Q}(H_{\Q}, \mathrm{gr}^2(\mathfrak{J}_H))$
be the homomorphism induced from the inclusion map $H_{\Q} \rightarrow \mathrm{gr}^1(\mathfrak{J}_H)$ defined by
\[ \sum_{l=1}^n c_{l} x_{l} \mapsto \sum_{l=1}^n c_{l} s_{11}(\overline{x_l}). \]
Let $p_2' : \mathrm{Hom}_{\Q}(H_{\Q}, \mathrm{gr}^2(\mathfrak{J}_H)) \rightarrow \mathrm{Hom}_{\Q}(H_{\Q}, S^2 H_{\Q})$
be the homomorphism induced from the projection $\mathrm{gr}^2(\mathfrak{J}_H) \rightarrow S^2 H_{\Q}$ defined by
\[ \sum_{t \in Y} c_t t \mapsto \sum_{1 \leq i \leq j \leq n} c_{\overline{v}_{11}(i,j)} x_i x_j. \]
Let $p_3' : \mathrm{Hom}_{\Q}(H_{\Q}, S^2 H_{\Q}) \rightarrow \mathrm{Hom}_{\Q}(H_{\Q}, H_{\Q}^{\otimes 2})$ be the homomorphism
induced from the homomorphism $S^2 H_{\Q} \rightarrow H_{\Q}^{\otimes 2}$ defined by
\[ x_i x_j \mapsto x_i \otimes x_j + x_j \otimes x_i \]
for any $1 \leq i \leq j \leq n$. Let $p_4' : \mathrm{Hom}_{\Q}(H_{\Q}, H_{\Q}^{\otimes 2}) \cong H_{\Q}^* \otimes_{\Q} H_{\Q}^{\otimes 2} \rightarrow H_{\Q}$
be the contraction map with respect to the first and the second component.
Then the composition map $p_4' \circ p_3' \circ p_2' \circ p_1'$, denoted by $p'$, is an $\mathrm{Aut}\,F_n$-equivariant homomorphism.
Set
\[ f_H := p' \circ \theta_H : \mathrm{Aut}\,F_n \rightarrow H_{\Q}. \]
Then we can see that
\[ f_H(\sigma) := \begin{cases}
                    0, \hspace{1em} & \sigma=P, Q, \\
                    -x_1, \hspace{1em} & \sigma=S, \\
                    -x_2, \hspace{1em} & \sigma=U.
                 \end{cases}\]
Hence, by the same argument as Theorem {\rmfamily \ref{T-Haruhi}}, we obtain the following.
\begin{thm}\label{T-Haruhi-H}
For any $n \geq 2$,
\[ f_M = - f_H + \delta_x \]
as a crossed homomorphism. Here $x := x_1 + x_2 + \cdots + x_n \in H_{\Q}$, and
$\delta_x$ is the principle derivation associated with $x$.
\end{thm}

We remark that we have
\[\begin{split}
  H^1(& \mathrm{Aut}\,F_n, \mathrm{Hom}_{\Q}(\mathrm{gr}^1(\mathfrak{J}_H), \mathrm{gr}^2(\mathfrak{J}_H))) \\
  & \cong H^1(\mathrm{Aut}\,F_n, H_{\Q}^* \otimes_{\Q} S^2 H_{\Q})^{\oplus \frac{1}{2}m^2(m^2-1)^2} \\
  & \hspace{4em} \bigoplus H^1(\mathrm{Aut}\,F_n, H_{\Q}^* \otimes_{\Q} \Lambda^2 H_{\Q})^{\oplus \frac{1}{2}(m^2-1)^2(m^2-4) \} } \\
  & \cong \Q^{\oplus \frac{1}{2}(m^2-1)^2(3m^2-8)}.
  \end{split}\]

\section{Acknowledgments}\label{S-Ack}

The part of this work was done when the author stayed at the Mathematical Institute of the University of Bonn as a visitor in 2017.
The author would like to thank the University of Bonn for its hospitality, and Tokyo University of Science for its finacial supports.
This work is supported by JSPS KAKENHI Grant Number 24740051.

\end{document}